%% file: neurips_2026.tex
\documentclass{article}

\usepackage[preprint]{neurips_2026}

\usepackage[utf8]{inputenc} 
\usepackage[T1]{fontenc}    
\usepackage{hyperref}       
\usepackage{url}            
\usepackage{booktabs}       
\usepackage{amsfonts}       
\usepackage{nicefrac}       
\usepackage{microtype}      
\usepackage{xcolor}         

\input{math_commands.tex}

\usepackage{amsmath}
\usepackage{amssymb}
\usepackage{mathtools}
\usepackage{amsthm}
\usepackage{algorithm}
\usepackage{algorithmic}
\usepackage{hyperref}
\usepackage{url}
\usepackage{amsmath, amsfonts, amsthm, amssymb, mathtools}
\usepackage{xspace}
\usepackage{enumitem}
\usepackage{mathrsfs, bbm}

\usepackage{comment}
\usepackage{threeparttable}
\usepackage{multirow}
\usepackage{multicol}
\usepackage{booktabs}
\usepackage{csquotes}
\usepackage{rotating}

\usepackage{graphicx}
\usepackage{wrapfig}
\usepackage{mathtools}

\theoremstyle{plain}
\newtheorem{theorem}{Theorem}[section]
\newtheorem{proposition}[theorem]{Proposition}
\newtheorem{lemma}[theorem]{Lemma}
\newtheorem{corollary}[theorem]{Corollary}
\theoremstyle{definition}
\newtheorem{definition}[theorem]{Definition}

\theoremstyle{remark}
\newtheorem{remark}[theorem]{Remark}
\newtheorem*{example}{Example}

\DeclarePairedDelimiter\ceil{\lceil}{\rceil}

\title{Smoothing Binary Optimization: A Primal-Dual Perspective}

\usepackage{authblk}
\author[1]{Wenbo Liu}
\author[2,3]{Akang Wang}
\author[1]{Dun Ma}
\author[4]{Hongyi Jiang}
\author[2]{Jianghua Wu}
\author[1]{Wenguo Yang}
\affil[1]{University of Chinese Academy of Sciences, China }
\affil[2]{Shenzhen Research Institute of Big Data, China}
\affil[3]{The Chinese University of Hong Kong, Shenzhen, China}
\affil[4]{City University of Hong Kong, Hong Kong SAR, China}

\begin{document}

\maketitle

\begin{abstract}
Binary optimization is a powerful tool for modeling combinatorial problems, yet scalable and theoretically sound solution methods remain elusive. Existing exact and heuristic solvers often face a trade-off between scalability and solution quality on large-scale instances. In this work, we introduce a novel primal-dual framework that reformulates unconstrained binary optimization as a continuous minimax problem, satisfying a strong max-min property. This reformulation effectively smooths the discrete problem, enabling the application of efficient gradient-based methods. We propose a simultaneous gradient descent-ascent algorithm that is highly parallelizable on GPUs and provably converges to a binary feasible solution in sublinear time. Extensive experiments on large-scale problems with up to 20,000 variables demonstrate that our method identifies high-quality solutions within seconds, significantly outperforming state-of-the-art alternatives.
\end{abstract}

\section{Introduction} 
\label{sec:Intro}

In this work, we consider the following \textit{binary optimization} problem:
\begin{equation}
\label{eq:binary_optimization}
\begin{aligned}
    &\underset{\vx\in\{0,1\}^n}{\min} \; F(\vx),
\end{aligned}
\end{equation}
where $\vx$ denotes a binary variable of dimension $n$ and $F:\{0,1\}^n\rightarrow\mathbb R$ represents an arbitrary real-valued function defined on the binary domain.
Problem~(\ref{eq:binary_optimization}) serves as a general model for a wide range of combinatorial optimization problems.
A prominent special case is the \textit{quadratic unconstrained binary optimization}~(QUBO), which is inherently equivalent to the Max-Cut problem and capable of characterizing all of the Karp's 21 $\mathcal {NP}$-complete problems~\citep{lucas2014ising}.
Moreover, \textit{integer linear programs} can be approximately transformed into QUBO by appropriately penalizing the constraints via quadratic terms~\citep{alessandroni2025alleviating}.
Other classic problems such as \textit{maximum satisfiability}~(MaxSAT)~\citep{bacchus2021maximum} are also encompassed within this framework.

Exact methods for binary optimization typically rely on branch-and-bound frameworks, integrated with various relaxation techniques such as the linear programming-based relaxations~\citep{barahona1989experiments} and the semi-definite relaxations~\citep{poljak1995solving}.
Recent advances have further enhanced these methods through sophisticated problem reduction techniques and accelerated cutting plane separation~\citep{rehfeldt2023faster}.
However, the inherent $\mathcal{NP}$-hardness fundamentally limits the scalability of exact approaches, which generally remain practical only for instances involving up to a few hundred variables.


Given the limitations of exact methods, research has shifted toward heuristics for large-scale binary optimization.
Traditional CPU-based methods involving customized tabu search~\citep{palubeckis2006iterated, wang2013probabilistic}, simulated annealing~\citep{bertsimas1993simulated} and genetic algorithms~\citep{merz1999genetic}, rely on well-crafted search mechanisms to iteratively explore the solution space.
However, their performance is often constrained by sequential execution.
Recent hardware progress has enabled significant improvements through GPU-acceleration, and conventional heuristics have benefited. 
For instance, ABS2~\citep{nakano2023diverse} integrates multiple search strategies and genetic operations executed in parallel on GPUs, while FEM~\citep{shen2025free} introduces the simulated annealing method from a continuous view and optimized via GPU implementations.
GPU acceleration has also propelled the development of first-order methods, yielding frameworks such as the sharp-peak penalized continuous relaxation by \cite{zhou2025sharp} and the sampling-based gradient policy by \cite{chen2025monte}.
Beyond classical computing approaches, quantum annealing on specialized hardware~\citep{berwald2019mathematics} has emerged as a promising alternative for QUBO.
However, the limited number of qubits and noise sensitivity~\citep{zaborniak2021benchmarking} still impede its widespread practical application.

Meanwhile, learning-based methods are gaining substantial traction as a powerful approach to binary optimization.
For example, QUBO problems can be naturally represented using graph structures, enabling the application of \textit{graph neural networks}~(GNNs) to effectively model variable interactions and achieve strong empirical performance
~\citep{schuetz2022combinatorial, ichikawa2024controlling}.
In addition, \cite{tonshoff2022one} introduces a general graph representation for \textit{constraint satisfaction problems} and proposes a GNN-based reinforcement learning approach enhanced with global search actions.
Despite their promise, such methods often operate without theoretical guarantees and encounter practical difficulties related to data acquisition and model generalization.

A prominent line of work within this domain \textit{smooths} Problem~(\ref{eq:binary_optimization}) using a probabilistic framework \citep{karalias2020erdos, qiu2024ros}: GNNs are leveraged to parameterize an independent distribution over the binary space $\{0,1\}^n$, and the loss function is subsequently constructed in terms of the expectation of $F$.
We note that this approach is essentially to relax Problem~(\ref{eq:binary_optimization}) to:
\begin{equation}
\label{eq:smooth_binary_optimization}
\min_{\vx \in [0,1]^n} f(\vx),
\end{equation}
where $f: \mathbb{R}^n \rightarrow \mathbb{R}$ is the \textit{multilinear extension} of $F$ (see Section~\ref{sec:preliminary} for definition). Although this relaxation preserves theoretical guarantees relative to the original problem~\citep{wang2022unsupervised}, the non-convex nature of $f$ in Problem~(\ref{eq:smooth_binary_optimization}) still presents significant optimization challenges.

\begin{figure*}[h]
\begin{center}
\includegraphics[width=0.9\textwidth]{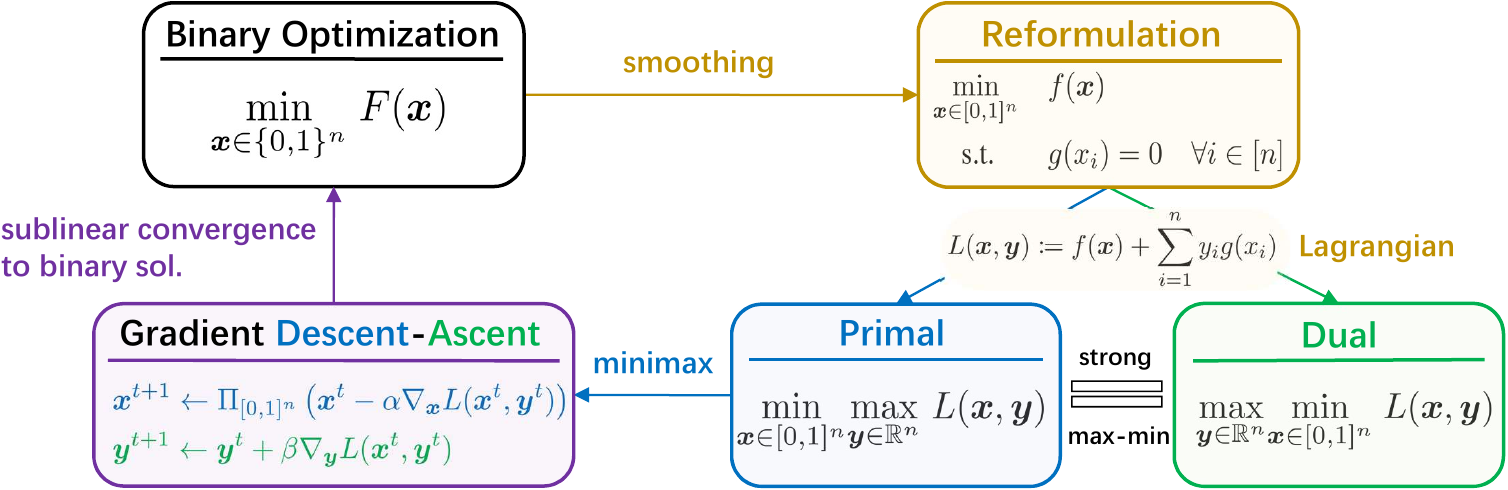}
\end{center}
\caption{An illustration of PDBO.}
\label{figure:PDBO}
\end{figure*}

To tackle this challenge, we revisit the binary optimization Problem~(\ref{eq:binary_optimization}), and reformulate it into a continuous model by characterizing binarity via continuous variables and convex constraint functions.
Based on Lagrangian duality, an equivalent \textit{minimax} problem is derived, to which we develop a \textit{gradient descent-ascent}~(GDA)-based algorithm for efficient solutions.
We remark that by initializing dual variables as positive values, we induce convexity in the primal space, promoting the identification of high-quality solutions.
As optimization proceeds, the dual variables gradually decrease, shifting emphasis from convexity promotion to the enforcement of integrality.
We refer to this approach as \enquote{a \underline{P}rimal-\underline{D}ual approach for \underline{B}inary \underline{O}ptimization}~(denoted as PDBO), with the overall framework illustrated in Figure~\ref{figure:PDBO}.

The distinct contributions of our work are summarized as follows.

\begin{itemize}[leftmargin=14pt]
    \item \textbf{Primal-Dual Reformulation.} We introduce a novel framework that reformulates unconstrained binary optimization as a continuous minimax problem, satisfying a strong max-min property. This reformulation smooths the discrete problem, enabling efficient gradient-based optimization.
    
     \item \textbf{GDA with Guarantees.} We develop a simultaneous GDA-based algorithm that is highly parallelizable on GPUs and converges to binary solutions in sublinear time.
    
    \item \textbf{Empirical Superiority.} Through extensive experiments on large-scale problems including Max-Cut, MaxSAT with up to 20,000 variables, we demonstrate that our method achieves state-of-the-art performance, obtaining high-quality solutions within seconds.
\end{itemize}

\section{Preliminaries} 
\label{sec:preliminary}
This section establishes foundational definitions—chief among them the \textit{multilinear extension}—and reviews its key properties for later use.
\begin{definition}[Extension]
\label{def:extension}
Let $F:\{0,1\}^n \rightarrow \mathbb{R}$ be a function defined on the binary domain. 
A function $f:[0,1]^n \rightarrow \mathbb{R}$ is called an \emph{extension} of $F$ if $f(\vx) = F(\vx)$ for every $\vx \in \{0,1\}^n$.
\end{definition}

\begin{definition}[Multilinear Extension]
    For any function $F:\{0,1\}^n \rightarrow \mathbb{R}$, its multilinear extension $f$ is defined by:
    \begin{equation*}
    f(\vx) \coloneqq \mathbb{E}_{\bm{\xi} \sim P_{\vx}}[F(\bm{\xi})] = \sum_{\bm{\xi} \in \{0,1\}^n} P_{\vx}(\bm{\xi}) \cdot F(\bm{\xi}),
\end{equation*}
where for $\vx \in [0,1]^n$, $P_{\vx}(\bm{\xi}) \coloneqq \prod_{i=1}^n x_i^{\xi_i}(1-x_i)^{1-\xi_i}$ denotes the $n$-dimensional independent Bernoulli distribution with mean vector $\vx$.
\end{definition}

We remark that $f$ is affine in each variable when all other variables are held fixed.
A key property of the multilinear extension is that it preserves the minimum value of the original function:

\begin{lemma}
\label{lem:multi_linear_rounding}
Let $f:[0,1]^n \rightarrow \mathbb{R}$ be the multilinear extension of $F:\{0,1\}^n \rightarrow \mathbb{R}$. Then
\begin{equation*}
    \min_{\vx \in [0,1]^n} f(\vx) = \min_{\vx \in \{0,1\}^n} F(\vx).
\end{equation*}
\end{lemma}

Lemma~\ref{lem:multi_linear_rounding} establishes the equivalence between Problem~(\ref{eq:binary_optimization}) and its continuous relaxation~(\ref{eq:smooth_binary_optimization}). However, we note that $f$ is generally non-convex, making the relaxed problem still challenging to solve.

\begin{proposition}
\label{prop:multi_linear_of_polynomial}
If $F:\{0,1\}^n \rightarrow \mathbb{R}$ is a polynomial function of the form $F(\vx) = \sum\limits_{\vd} a_{\vd} \vx^{\vd}$, where $\vx^{\vd} \coloneqq x_1^{d_1} x_2^{d_2} \cdots x_n^{d_n}$, then its multilinear extension is given by:
\begin{equation*}
    f(\vx) = \sum_{\vd} a_{\vd}  \prod_{i: d_i \neq 0} x_i.
\end{equation*}
\end{proposition}

Proposition~\ref{prop:multi_linear_of_polynomial} shows that for polynomial objectives (such as QUBO and MaxSAT), the multilinear extension has a compact closed form, avoiding exponential summation. This makes computation tractable for many practical optimization problems.

\begin{definition}[Minimax Problem and Critical Points]
\label{def:minimax}
Consider the minimax problem
\begin{equation*}
    \min_{\vx \in \mathcal{X}} \max_{\vy \in \mathcal{Y}} L(\vx, \vy),
\end{equation*}
where $\mathcal{X} \subseteq \mathbb{R}^n$ and $\mathcal{Y} \subseteq \mathbb{R}^m$ are the domains of $\vx$ and $\vy$, respectively. For a point $(\vx^*, \vy^*) \in \mathcal{X} \times \mathcal{Y}$:

\begin{itemize}[leftmargin=10pt]
    \item $(\vx^*, \vy^*)$ is a \emph{saddle point} if for all $(\vx, \vy) \in \mathcal{X} \times \mathcal{Y}$:
        $L(\vx^*, \vy) \leq L(\vx^*, \vy^*) \leq L(\vx, \vy^*)$.
    
    \item $(\vx^*, \vy^*)$ is a \emph{stationary point} if $\nabla_\vx L(\vx^*, \vy^*) = 0$ and $\nabla_\vy L(\vx^*, \vy^*) = 0$.
\end{itemize}
\end{definition}

\begin{remark}
\label{rem:saddle_existence}
A saddle point does not necessarily exist for an arbitrary minimax problem. However, if the strong max-min property holds, i.e.,
\begin{equation*}
    \min_{\vx \in \mathcal{X}} \max_{\vy \in \mathcal{Y}} L(\vx, \vy) = \max_{\vy \in \mathcal{Y}} \min_{\vx \in \mathcal{X}} L(\vx, \vy) \in \mathbb{R},
\end{equation*}
then there exists at least one saddle point.
\end{remark}

\section{Methodology}
\label{sec:method}

In Section~\ref{subsec:primal_dual_framework}, we introduce a primal-dual framework that reformulates the binary optimization Problem~(\ref{eq:binary_optimization}) as a minimax problem. Then, in Section~\ref{subsec:minimax_optimization}, we develop a GDA-based algorithm to solve the resulting minimax formulation and establish its convergence guarantees.

\subsection{Primal-Dual Reformulation}
\label{subsec:primal_dual_framework}

We begin by reformulating Problem~(\ref{eq:binary_optimization}) as a constrained continuous optimization problem:
\begin{equation}
\label{eq:qubo_conti}
\begin{aligned}
&\min_{\vx \in [0,1]^n} && f(\vx) \\
&\quad \text{s.t.} && g(x_i) = 0 && \forall i \in [n],
\end{aligned}
\end{equation}
where $f$ is the multilinear extension of $F$, and $g: [0,1] \rightarrow \mathbb{R}$ is a  constraint function to enforce binarity such that the condition $g(x_i)=0$ is equivalent to $x_i\in\{0,1\}$.

An appropriate function $g$ should satisfy the following sufficient conditions:
\begin{itemize}
    \item $g$ is convex and continuous on $[0,1]$, with $g(0) = g(1) = 0$;
    \item $g$ is differentiable on $(0,1)$, and $g'(x) = 0$ if and only if $x = \frac{1}{2}$;
    \item $g$ is symmetric on $[0,1]$, i.e., $g(x) = g(1-x)$ for all $x \in [0,1]$.
\end{itemize}

These conditions imply that $g(x)< 0$ for all $x\in(0,1)$ and $g'(x)< 0$ for all $x\in(0,\frac12)$, and therefore $-g$ can be regarded as a metric of fractionality.
Several choices of $g$ have been proposed in the literature. Common examples include an even-degree polynomial $g(x)\coloneqq (2x-1)^{2d} - 1$~\citep{ichikawa2024controlling} and the entropy-based function $g(x)\coloneqq x\log(x)+(1-x)\log(1-x)$~\citep{shen2025free}, extended continuously so that $g(0)=g(1)=0$.

Let $\vy \in \mathbb{R}^n$ denote the dual variables associated with the constraints in Problem~(\ref{eq:qubo_conti}). We define the Lagrangian function as:
\begin{equation*}
\label{eq:lagrangian}
L(\vx,\vy) \coloneqq f(\vx) + \sum_{i=1}^n y_i g(x_i).
\end{equation*}

\begin{proposition} \label{prop:lagrangian_primal}
Let $v^*$ denote the optimal value of Problem~(\ref{eq:qubo_conti}), then
$v^* = \underset{\vx\in[0,1]^n}{\inf}\underset{\vy\in\mathbb R^n}{\sup}\;L(\vx,\vy)$.
\end{proposition}

Proposition~\ref{prop:lagrangian_primal} characterizes the primal problem as a minimax optimization. The following theorem establishes strong duality for Problem~(\ref{eq:qubo_conti}).

\begin{theorem}
\label{thm:strong_duality}
\begin{equation}
\label{eq:strong_duality_inf}
\inf_{\vx \in [0,1]^n} \sup_{\vy \in \mathbb{R}^n} L(\vx,\vy) = \sup_{\vy \in \mathbb{R}^n} \inf_{\vx \in [0,1]^n} L(\vx,\vy).
\end{equation}
\end{theorem}

Equation~(\ref{eq:strong_duality_inf}) in Theorem~\ref{thm:strong_duality}, known as the \textit{strong max-min property}~\citep{boyd2004convex}, implies the existence of a \textit{saddle point} of the Lagrangian function that corresponds to an optimal solution of Problem~(\ref{eq:qubo_conti}).
We therefore reformulate Problem~(\ref{eq:qubo_conti}) as the following equivalent minimax problem:
\begin{equation}
\label{eq:minimax_problem}
\min_{\vx \in [0,1]^n} \max_{\vy \in \mathbb{R}^n} L(\vx,\vy),
\end{equation}
where the order of minimization and maximization can be interchanged.

\subsection{Minimax Optimization}
\label{subsec:minimax_optimization}
The Lagrangian $L(\vx,\vy)$ is linear in $\vy$ but not necessarily convex in $\vx$. Consequently, solving the minimax problem~(\ref{eq:minimax_problem}) to global optimality remains $\mathcal{NP}$-hard.
In this section, we develop a GDA-based algorithm for solving the minimax problem~(\ref{eq:minimax_problem}) and establish its convergence guarantees.

\subsubsection{Gradient Descent-Ascent}
\label{subsec:gda}

For Problem~(\ref{eq:minimax_problem}), a classical approach is the GDA method, which alternates between gradient descent on the primal variable $\vx$ and gradient ascent on the dual variable $\vy$. 
The strong max-min property enables \textit{simultaneous} updates as follows:
\begin{equation}
\label{eq:GDA_for_QUBOdd}
\begin{aligned}
&\vx^{t+1} \gets \Pi_{[0,1]^n}\left(\vx^t - \alpha \cdot \nabla_{\vx} L(\vx^t, \vy^t) \right), \\
&\vy^{t+1} \gets \vy^t + \beta \cdot \nabla_{\vy} L(\vx^t, \vy^t),
\end{aligned}
\end{equation}
where $\vx^t$ and $\vy^t$ denote the primal and dual variables at iteration $t$, $\Pi_{[0,1]^n}$ denotes the projection onto the hypercube $[0,1]^n$, and $\alpha, \beta > 0$ are the step sizes.

Notably, the projection ensures $\vx^t \in [0,1]^n$ throughout the optimization, which implies the non-increasing monotonicity of $\{y_i^t\}_t$ since $g(x_i^t) \leq 0$. Moreover, the term $\sum_{i=1}^n (-y_i) \cdot (-g(x_i))$ in $L(\vx,\vy)$ can be interpreted as an integrality penalty, where the penalty coefficients $-y_i$ increase monotonically during optimization.
Moreover, GDA introduces a principled update rule for the penalty coefficients, distinguishing PDBO from traditional penalty methods.
However, increasing penalties alone does not guarantee integer solutions, as fractional stationary points may persist.

\begin{example}[Max-Cut]
Let $f(\vx) \coloneqq \vx^\top \mW \vx - \vone^\top \mW \vx$ and $g(x_i) \coloneqq x_i^2 - x_i$. Then the fractional solution $\vx^t \coloneqq (\frac{1}{2}, \ldots, \frac{1}{2}) \in \mathbb{R}^n$ yields $\nabla_\vx L(\vx^t,\vy^t) = 0$ for any $\vy^t \in \mathbb{R}^n$.
\end{example}

To overcome this limitation, we incorporate targeted perturbations to prevent convergence to fractional points. Specifically, given an initial solution $(\vx^0,\vy^0) \in [0,1]^n \times \mathbb{R}^n_{++}$, maximum iterations $T_{\max}$, and tolerance $\delta > 0$, we iteratively perform the updates in~(\ref{eq:GDA_for_QUBOdd}). When $x_i^t$ is closed to and tends to stagnate near $\frac12$, we perturb it to maintain a minimum distance of $\delta$ from $\frac{1}{2}$. The complete procedure is described in Algorithm~\ref{alg:PDBO}.

\begin{algorithm}[htbp]
\caption{PDBO}
\label{alg:PDBO}
\begin{algorithmic}[1]
\REQUIRE Initial solution  $ (x^0, y^0) \in [0,1]^n \times \mathbb{R}^n_{++} $ , 
         Stepsizes  $ (\alpha, \beta) \in \mathbb{R}_{++}^2 $ , 
         Tolerance  $ \delta > 0 $ , 
         Number of iterations  $ T_{\text{max}} $ 

\FOR{ $ t = 0, 1, 2, \dots, T_{\text{max}} - 1 $ }
    \FOR{ $ i = 1, 2, \dots, n $ }
        \IF{ $ \left|x_i^t - \frac{1}{2}\right| \leq \delta $  \AND  $ \left|\frac{\partial}{\partial x_i} L(\mathbf{x}^t, \mathbf{y}^t)\right| \leq 2\delta $  \AND  $ y_i^t \leq 0 $ }
            \STATE  $ 
                x_i^{t+1} \gets 
                \begin{cases}
                    \frac{1}{2} - \delta & \text{if } x_i^t \leq \frac{1}{2}, \\
                    \frac{1}{2} + \delta & \text{otherwise}.
                \end{cases}
             $ 
        \ELSE
            \STATE  $ x_i^{t+1} \gets \Pi_{[0,1]}\left( x_i^t - \alpha \cdot \frac{\partial}{\partial x_i} L(\mathbf{x}^t, \mathbf{y}^t) \right) $ 
        \ENDIF
        \STATE  $ y_i^{t+1} \gets y_i^t + \beta \cdot g(x_i^t) $ 
    \ENDFOR
\ENDFOR
\end{algorithmic}
\end{algorithm}


\subsubsection{Convergence Analysis}
\label{subsubsec:convergence_analysis}
This section analyzes the convergence of the primal variable $\vx$ in Algorithm~\ref{alg:PDBO} to a binary point, and provides corresponding complexity guarantees.
We begin by establishing a lower bound for $\evy_i^t$.
\begin{proposition}
\label{prop:y_bounded}
Define $\Theta \coloneqq \underset{\vx\in[0,1]^n}\max||\nabla_\vx f(\vx)||_1$ and  $b \coloneqq \frac{1+\Theta}{g'(\frac12-\delta)} + \left(2+ \ceil{\frac{1}{2\alpha}}  \right)\cdot \beta\cdot g(\frac12)$.
Then for each $i\in[n]$ and any $t \geq 0$, we have $\evy_i^t\geq b$.
\end{proposition}
Note that for any $i\in[n]$, the sequence $\{y_i^t\}_{t\geq 0}$ is monotone non-increasing, and hence the lower bound leads to the convergence of $\{\evy_i^t\}_{t\geq 0}$.
Furthermore, $\{g(\evx_i^t)\}_{t\geq 0}$ converges to zero since $\beta g(\evx_i^t)=\evy_i^{t+1}-\evy_i^t$.
\begin{corollary}
    The sequence $\{\vy^t\}_{t\geq 0}$ converges to some $\vy^*\in\mathbb R^n$ with $\evy_i^*\geq b$ for any $i\in[n]$.
\end{corollary}

\begin{corollary}
\label{cor:x_convergence}
    For any $i\in[n]$, the sequence $\{g(\evx_i^t)\}_{t\geq 0}$ converges to $0$.
\end{corollary}
We remark that Corollary~\ref{cor:x_convergence} implies the convergence of $\vx$ to a binary point.
This holds because $||\nabla_\vx L(\vx^t,\vy^t)||$ is bounded, and by choosing a sufficiently small step size $\alpha$, oscillations between binary points can be avoided.

We now establish the iteration complexity for Algorithm~\ref{alg:PDBO} to reach a binary solution.
Since the primary goal is to solve the binary optimization problem~(\ref{eq:binary_optimization}), we focus on convergence toward integrality rather than traditional stationarity criteria.
The concept is formalized as follows:
\begin{definition}
Given Problem~(\ref{eq:qubo_conti}), a point $\vx\in[0,1]^n$ is called an $\epsilon$-binary point if $-\sum\limits_{i=1}^ng(x_i) \leq \epsilon$.
\end{definition}

\begin{theorem}
\label{thm:convergence_rate}
    The iteration complexity for Algorithm~\ref{alg:PDBO} to return an $\epsilon$-binary point is bounded by 
    \begin{equation*}
        \mathcal{O}\left(\frac{\Vert\vy^0-\vy^*\Vert_1}{\beta\epsilon}\right).
    \end{equation*}
\end{theorem}
Theorem~\ref{thm:convergence_rate} shows that Algorithm~\ref{alg:PDBO} converges to a binary feasible point at a sublinear rate.
While the convergence analysis establishes binarity rather than objective optimality, Section~\ref{sec:landscape_analysis} complements it with a landscape analysis that helps explain the superior empirical performance of PDBO.

\section{Numerical Experiments}
\label{sec:experiments}
We conduct numerical experiments to evaluate the performance of our proposed \texttt{PDBO} method. We evaluate \texttt{PDBO} on three application domains that fit our problem framework, as summarized in Section~\ref{subsection:application}.
Our code is available at: \url{https://anonymous.4open.science/r/PD-Bo/}.

\subsection{Applications}
\label{subsection:application}
We introduce three classic combinatorial optimization problems that align with Problem~(\ref{eq:binary_optimization}): Max-Cut, Max-$k$-SAT, and Max-$k$-Cut.
Additional experiments on Maximum Independent Set (MIS) and Maximum a Posteriori (MAP) inference are presented in Appendix~\ref{apdx:additional_experiments}.

\textbf{Max-Cut.} Given a graph $G\coloneqq(V,E)$ with adjacency matrix $\mathbf{W}$, the Max-Cut problem seeks to partition the vertex set $V$ into two disjoint subsets that maximize the total weight of edges between them. This can be formulated as the following binary optimization problem:
\vskip-0.1in
\begin{equation*}
    \max_{\vx \in \{0,1\}^n} \sum_{(i,j) \in E} W_{ij} (x_i + x_j - 2x_i x_j),
\end{equation*}
\vskip-0.1in
where $x_i \in \{0,1\}$ indicates the partition assignment of vertex $i$, and $W_{ij}$ denotes the weight of edge $(i,j) \in E$. The expression $x_i + x_j - 2x_i x_j$ equals 1 if $x_i \neq x_j$ and 0 otherwise.


\textbf{Max-\textit{k}-SAT.} 
Given a conjunctive normal form (CNF) formula with $m$ clauses, where each clause contains at most $k$ literals, the goal is to find a truth assignment that maximizes the number of satisfied clauses. This problem can be formulated as the following binary optimization problem:
\vskip-0.2in
\begin{equation*}
\begin{aligned}
    &\min_{\vx \in \{0,1\}^n} \; \sum_{j=1}^{m} \prod_{l_i \in C_j} p(x_i)
\end{aligned}
\end{equation*}
\vskip-0.1in

where $x_i \in \{0,1\}$ denotes the truth assignment of the $i$-th variable, $C_j$ denotes the $j$-th clause, and $p(x_i) = x_i$ if literal $l_i$ appears negated in clause $C_j$ and $1- x_i$ otherwise.
The product term $\prod_{l_i \in C_j} p(x_i)$ equals 1 if clause $C_j$ is unsatisfied and 0 if it is satisfied. 

\textbf{Max-\textit{k}-Cut.} 
The Max-\textit{k}-Cut problem generalizes Max-Cut by partitioning the vertex set $V$ into $k$ disjoint subsets, maximizing the total weight of edges crossing between different subsets. The discrete formulation is:
\vspace{-0.1cm}
\begin{equation*}
\label{eq:max_k_cut_discrete}
\begin{aligned}
&\max_{\mX \in \mathbb{R}^{k \times n}} && \sum_{(i,j) \in E} W_{ij} \left(1 - \mX_{:,i}^\top \mX_{:,j}\right) \\
&\quad \text{s.t.} && \mX_{:,i} \in \{\ve^{(1)}, \ldots, \ve^{(k)}\}, && \forall i \in [n]
\end{aligned}
\end{equation*}
where $\mX_{:,i}$ is a $k$-dimensional one-hot vector indicating the partition assignment of vertex $i$ since $\ve^{(j)}$ denotes the $j$-th standard basis vector. The expression $1 - \mX_{:,i}^\top \mX_{:,j}$ equals 1 if vertices $i$ and $j$ belong to different subsets and 0 otherwise.
We remark that Max-\textit{k}-Cut can be incorporated into our primal-dual framework with a natural generalization. Please see Appendix~\ref{apdx:maxkcut} for technical details.

\subsection{Setup}
\label{subsec:experimental_setup}

\textbf{Datasets.}
For Max-Cut and Max-\textit{k}-Cut, we use the Gset benchmark~\citep{ye2003}, a widely recognized dataset comprising graphs with $800$ to $20{,}000$ nodes, totaling 71 instances.
For Max-\textit{k}-SAT, we use the datasets from~\cite{tonshoff2022one}, which include 3CNF, 4CNF, and 5CNF formulas, each containing 50 instances. This enables consistent and comparable evaluation with existing work.

\textbf{Baselines.}
We consider two types of baseline algorithms: For traditional heuristics with GPU accelerations, we adopt (i) \texttt{ABS2}~\citep{nakano2023diverse}: an evolutionary method featuring three diverse strategies: multiple search, multiple genetic operations and multiple solution pools; and (ii) \texttt{FEM}~\citep{shen2025free}: an annealing-based method with additional techniques tailored for Max-Cut.
For the learning-based approaches, we consider the following state-of-the-art methods: (iii) \texttt{PIGNN}~\citep{schuetz2022combinatorial}: an unsupervised method for QUBO problems, which leverage physics-inspired GNNs and delivers commendable performance; (iv) \texttt{ANYCSP}~\citep{tonshoff2022one}: a GNN-based reinforcement learning approach for constraint satisfaction problems, utilizing a compact graph representation and global search actions; (v) \texttt{CRA}~\citep{ichikawa2024controlling}: an annealing-based method that control the convexity of objectives, with GNNs leveraged to further enhance the performance; and (vi) \texttt{ROS}~\citep{qiu2024ros}: a GNN-guided relax-optimize-and-sample method with fine tuning techniques, designed for Max-$k$-Cut problems. Finally, we consider the advanced commercial solver (vii) \texttt{Gurobi}~12.0.2~\citep{gurobi}.

\textbf{Model Settings.}
To enforce binarity in Problem~(\ref{eq:qubo_conti}), we employ the function $g(x) \coloneqq x^2 - x$ for \texttt{PDBO} unless otherwise specified. To fully utilize GPU parallelization, we sample $B$ independent initial solutions and execute \texttt{PDBO} in parallel from each starting point, retaining the best solution across all runs. Detailed parameter configurations, including the values of $B$, $\vy^0$, $\alpha$, and $\beta$, are provided in Appendix~\ref{apdx:hyper_parameters}, together with a sensitivity analysis of the key hyperparameters.

\textbf{Evaluation Configuration.}
All experiments are conducted with a 180-second time limit using identical hardware configuration: a 12th Gen Intel Core i9-12900K CPU and an NVIDIA GeForce RTX 3090 GPU. 
\textbf{For each method, we report the solution time as the earliest time at which it attains its best solution within the 180-second budget.}
\texttt{PDBO} is implemented in JAX 0.6.1~\citep{jax}.

\subsection{Computational Results}
\subsubsection{Results for Max-Cut}

We evaluate \texttt{PDBO} by measuring both solution quality (number of cut edges) and computational efficiency (runtime in seconds). 
Our main results focus on the five largest Gset graphs (G67, G70, G72, G77, and G81), each containing over 10,000 nodes, which present significant computational challenges for state-of-the-art methods. Results for smaller graphs are provided in Appendix~\ref{apdx:additional_experiments}.

\begin{table}[h]
\caption{Results on Gset instances for Max-Cut}
\label{tab:maxcut}
\centering
\resizebox{\textwidth}{!}{
\begin{tabular}{@{}lcccccccccccc@{}}
\toprule
\multirow{2}{*}{Method} & \multicolumn{2}{c}{G67 (n=10k)} & \multicolumn{2}{c}{G70 (n=10k)} & \multicolumn{2}{c}{G72 (n=10k)} & \multicolumn{2}{c}{G77 (n=14k)} & \multicolumn{2}{c}{G81 (n=20k)} \\
\cmidrule(lr){2-3} \cmidrule(lr){4-5} \cmidrule(lr){6-7} \cmidrule(lr){8-9} \cmidrule(lr){10-11}
 & Obj $\uparrow$ & Time $\downarrow$ & Obj $\uparrow$ & Time $\downarrow$ & Obj $\uparrow$ & Time $\downarrow$ & Obj $\uparrow$ & Time $\downarrow$ & Obj $\uparrow$ & Time $\downarrow$\\
\midrule[0.4 pt]
\texttt{PIGNN} &  5538 & 23.7 & 8534 & 25.2 & 5588 & 44.5 & 7896 & 42.1 & 11078 & 157.6 \\
\texttt{CRA} &  5948 & 53.7 & 9240 & 51.7 & 6058 & 53.9 & 8720 & 75.9 & 12450 & 120.4 \\
\texttt{ROS} &  6144 & 1.5 & 8872 & 1.8 & 6148 & 1.2 & 8746 & 2.2 & 12320 & 5.2\\
\texttt{ANYCSP} &  6772 & 39.9 & 9379 & 35.7 & 6816 & 36.1 & 9686 & 53.5 & 13670 & 73.5 \\
\texttt{FEM} &  6782 & 2.4 & 5120 & 0.2 & 6824 & 2.6 & 9688 & 4.0 & 13684 & 7.5 \\
\texttt{ABS2} &  6880 & 156.3 & 9510 & 175.5 & 6932 & 172.2 & 9824 & 171.1 & 13850 & 177.1 \\
{\texttt{Gurobi}}&{\bf{6944}}&{51.9}&{9514}&{133.3}&{\bf{6990}}&{171.6}&{\bf{9882}}&{175.1}&{13848}&{179.8}\\
\texttt{PDBO} &  6872 & 2.5 & \textbf{9537} & 1.9 & 6906 & 2.0 & 9812 & 2.1 & \textbf{13852} & 3.3 
 \\
\bottomrule
\end{tabular}}
\end{table}

Table~\ref{tab:maxcut} summarizes the results for all baselines. The \enquote{Obj.} column shows the best objective value found by each method within the 180-second limit. The \enquote{Time} entry corresponds to the earliest time at which the best solution was found.
Note that $\uparrow$ denotes that larger values are better, and $\downarrow$ indicates that smaller values are preferred.
The results show that \texttt{PDBO} exhibits substantially greater efficiency than other baselines, and delivers the best solutions for G70 and G81.
Meanwhile, \texttt{Gurobi} finds the best solutions to the other three instances. It is worth noting, as highlighted in~\cite{rehfeldt2023faster}, that recent versions of Gurobi have seen major improvements on QUBO and Max-Cut problems, likely due to the integration of specialized heuristics.

\begin{wrapfigure}{r}{0.45\textwidth}
\vskip -0.3 in
    \centering
    \includegraphics[width=0.9\linewidth]{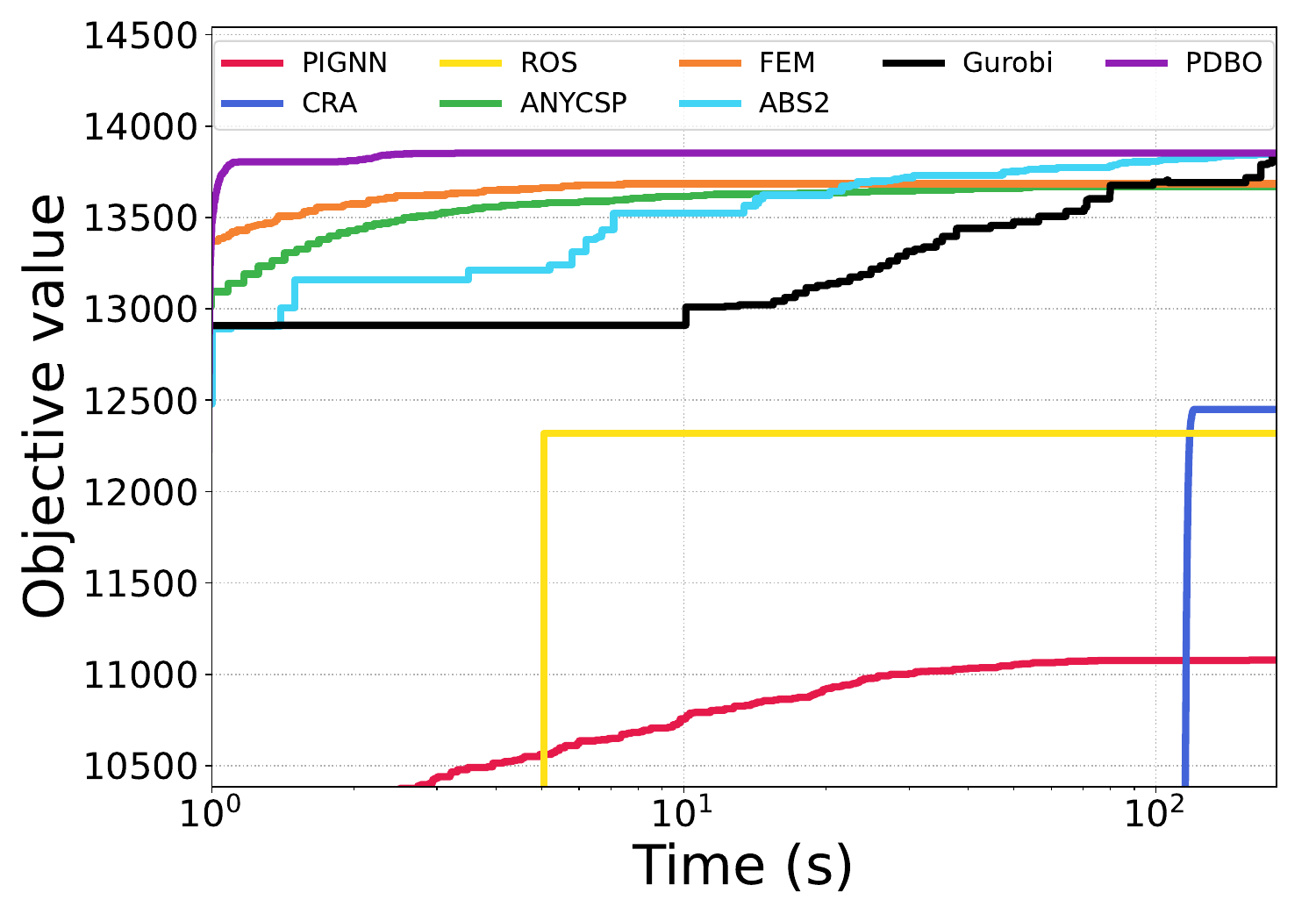}
    \vskip-0.1in
    \caption{The objective value versus runtime for G81}
    \label{fig:maxcut}
    \vskip-0.15 in
\end{wrapfigure}

Figure~\ref{fig:maxcut} compares the solution quality throughout the optimization process, plotting the objective value against runtime for two representative instances.
Clearly, \texttt{PDBO} (purple) achieves and maintains the highest objective values shortly after initialization, demonstrating both rapid convergence and superior solution quality for most of the time horizon.
Although \texttt{Gurobi} (black) eventually attains a marginally better solution given the full time budget, \texttt{PDBO} maintains a substantial advantage throughout the majority of the optimization timeline.



\subsubsection{Results for Max-$k$-SAT.}

\begin{wraptable}[7]{r}{8.2cm}
\vskip -0.5cm
    \caption{Results on CNF instances for Max-$k$-SAT}
    \vskip 0.1 in
    \label{tab:maxsat}
    \centering
    \resizebox{\linewidth}{!}{
    \begin{tabular}{@{}lcccccccc@{}}
\toprule
\multirow{2}{*}{Method} & \multicolumn{2}{c}{3CNF} & \multicolumn{2}{c}{4CNF} & \multicolumn{2}{c}{5CNF} \\
\cmidrule(lr){2-3} \cmidrule(lr){4-5} \cmidrule(lr){6-7} 
& Obj $\downarrow$ & Time $\downarrow$ & Obj $\downarrow$ & Time $\downarrow$ & Obj $\downarrow$ & Time $\downarrow$  \\
\midrule
\texttt{FEM} & 1885.2$_{\pm23.9}$ & 0.8$_{\pm0.1}$ & -- & -- &--&--\\
\texttt{ANYCSP} &  1583.3 $_{\pm 17.5 }$ & 123.8 $_{\pm 40.7 }$ & 1210.9 $_{\pm 12.6 }$ & 141.6 $_{\pm 25.9 }$ & 1213.7 $_{\pm 11.4 }$ & 141.0 $_{\pm 31.9 }$ & \\
{\texttt{Gurobi}} &{9322.6 $_{\pm  64.2 }$} &{0.0 $_{\pm  0.1 }$
}&{9329.2 $_{\pm  69.7 }$ }&{0.0 $_{\pm  0.0 }$ }&{9310.7 $_{\pm  75.4 }$ }&{0.1 $_{\pm  0.0 }$ }\\
\texttt{PDBO} &  \textbf{1582.7} $_{\pm 12.0 }$ & 0.9 $_{\pm 0.1 }$ & \textbf{1147.4} $_{\pm 0.8 }$ & 1.2 $_{\pm 0.1 }$ & \textbf{976.0} $_{\pm 9.3 }$ & 2.2 $_{\pm 0.3}$ & \\
\bottomrule
\end{tabular}}
\vskip 0.1 in
\end{wraptable}

We evaluate \texttt{PDBO} against baseline algorithms \texttt{FEM}, \texttt{ANYCSP} and \texttt{Gurobi} on Max-\textit{k}-SAT using the 3CNF, 4CNF, and 5CNF datasets, each containing 50 problem instances. 
Since the objective in Max-$k$-SAT is a polynomial of higher degree, we accordingly adopt $g(x) \coloneqq x \log x + (1-x) \log(1-x)$ in \texttt{PDBO} to enforce binarity, which yields better performance. 
Performance is measured by the average number of unsatisfied clauses across all instances in each dataset, as summarized in Table~\ref{tab:maxsat}, where entries marked with \enquote{--} indicate memory issues.
We find that both \texttt{FEM} and \texttt{Gurobi} struggle on these problems, yielding inferior solutions compared to \texttt{ANYCSP} and \texttt{PDBO}.
Moreover, the results show that \texttt{PDBO} achieves competitive or superior performance compared to \texttt{ANYCSP}, while avoiding the substantial training overhead and significantly longer inference time required by the latter.

\subsubsection{Results for Max-$k$-cut}

We evaluate the performance of \texttt{PDBO} against baseline algorithms \texttt{ROS}, \texttt{ANYCSP}, \texttt{FEM} and \texttt{Gurobi} on Gset instances for the Max-3-Cut problem. Following our evaluation protocol, we present results for instances with more than 10,000 nodes in Table~\ref{tab:max3cut}, with remaining results provided in Appendix~\ref{apdx:additional_experiments}.
Table~\ref{tab:max3cut} demonstrates that \texttt{PDBO} consistently outperforms all baseline methods across all evaluated Max-3-Cut instances, achieving superior solution quality with competitive runtime efficiency.
\begin{table}[h]
    \caption{Results on Gset instances for Max-3-Cut}
    \label{tab:max3cut}
    \centering
    \resizebox{\textwidth}{!}{
\begin{tabular}{@{}lcccccccccccc@{}}
\toprule[0.5pt]
\multirow{2}{*}{Method} & \multicolumn{2}{c}{G67 (n=10k)} & \multicolumn{2}{c}{G70 (n=10k)} & \multicolumn{2}{c}{G72 (n=10k)} & \multicolumn{2}{c}{G77 (n=14k)} & \multicolumn{2}{c}{G81 (n=20k)} \\
\cmidrule(lr){2-3} \cmidrule(lr){4-5} \cmidrule(lr){6-7} \cmidrule(lr){8-9} \cmidrule(lr){10-11}
 & Obj $\uparrow$ & Time $\downarrow$ & Obj $\uparrow$ & Time $\downarrow$ & Obj $\uparrow$ & Time $\downarrow$ & Obj $\uparrow$ & Time $\downarrow$ & Obj $\uparrow$ & Time $\downarrow$ \\
\midrule[0.4pt]
\texttt{ROS} &  7364 & 3.0 & 9983 & 1.9 & 7435 & 2.9 & 10559 & 5.3 & 14907 & 9.7 \\
\texttt{ANYCSP} &  7797 & 55.9 & 9909 & 16.2 & 7906 & 58.9 & 11158 & 84.0 & 15727 & 115.0 \\
\texttt{FEM} &  7748 & 3.3 & \textbf{9999} & 1.4 & 7835 & 0.7 & 11102 & 4.4 & 15683 & 6.1 \\
{\texttt{Gurobi}} & {7855}&{178.1}& {\textbf{9999}}&{0.7}&{8045}&{174.6}&{11101}&{179.9}&{15146}&{179.9}\\
\texttt{PDBO} &  \textbf{8015} & 6.0 & \textbf{9999} & 3.3 & \textbf{8111} & 4.7 & \textbf{11467} & 5.0 & \textbf{16191} & 7.9 \\
\bottomrule[0.5pt]
\end{tabular}}
\end{table}

\subsection{Landscape Analysis}
\label{sec:landscape_analysis}

We next use Max-Cut to illustrate how the primal-dual dynamics shape the optimization landscape.
Recall the minimization form of Max-Cut, $f(\vx)=\vx^\top \mW \vx-\vone^\top \mW \vx$, 
where $\mW$ is symmetric with zero diagonal entries, so the multilinear extension coincides with the quadratic objective itself.
With $g(x_i)=x_i^2-x_i$, the associated Lagrangian is
    $L(\vx,\vy)
    =\vx^\top \mW \vx-\vone^\top \mW \vx
    +\sum_{i=1}^n y_i(x_i^2-x_i)$
and its Hessian with respect to $\vx$ is
\begin{equation*}
    \nabla_{\vx}^2 L(\vx,\vy)=2\big(\mW+\operatorname{diag}(\vy)\big).
\end{equation*}




This simple form reveals a convex-to-nonconvex landscape evolution. 

\textit{Convex initialization:} 
PDBO initializes the dual variables as $\vy^0=\bar y\vone$; under the condition $\bar y\geq -\lambda_{\min}(\mW)$, one has $\mW+\operatorname{diag}(\vy^0)\succeq 0$, which renders $L(\cdot,\vy^0)$ globally convex. 
The method therefore starts from a single-basin landscape without spurious local minima. 

\textit{Gradual nonconvexification:} 
The dual update obeys $y_i^{t+1}-y_i^t=\beta\big((x_i^t)^2-x_i^t\big)\leq 0$, so each coordinate of $\vy^t$ decreases monotonically. 
As a result, the diagonal curvature in $\operatorname{diag}(\vy^t)$ is progressively weakened, and the landscape shifts from convex to nonconvex as PDBO proceeds. 

\textit{Binary recovery:} 
The dual term can be rewritten as $y_i(x_i^2-x_i)=-y_i(x_i-x_i^2)$, where $x_i-x_i^2\geq 0$ measures the fractionality of $x_i$. 
The decrease of $y_i$ increases the coefficient $-y_i$, thereby inducing a stronger penalty against fractional values and driving the iterates toward binary solutions, in agreement with the convergence analysis provided in Section~\ref{subsubsec:convergence_analysis}.


\textbf{Empirical evidence.}
Table~\ref{tab:landscape_ablation} further supports this interpretation through a sensitivity analysis on the initialization magnitude $\bar y$.
Each entry reports the obtained Max-Cut objective value.
Here, \enquote{$\surd$} indicates that $\mW+\bar y \cdot \mI \succeq 0$, while \enquote{$\times$} indicates otherwise.
The results show that the objective value improves substantially once $\bar y$ is large enough to make the initial landscape convex.

\begin{table}[h]
\caption{Effect of the initial dual magnitude $\bar y$ on Max-Cut objective values.}
\label{tab:landscape_ablation}
\centering
\resizebox{\linewidth}{!}{
\begin{tabular}{lccccccc}
\toprule[0.5pt]
Instance & $\bar y=0 (\times)$ & $\bar y=1 (\times)$ & $\bar y=2 (\times)$ & $\bar y=4 (\surd)$ & $\bar y=6 (\surd)$ &$\bar y=8 (\surd)$ & $\bar y=10 (\surd)$ \\
\midrule[0.4pt]
G67 & 5632  & 6394 &6674  &6864  &{6872}  &6868   &6868  \\
G70 & 8756  & 9145  & 9398  & 9535  & {9537}  & 9536  & 9536  \\
G72 &5714  &6474   &6720  & 6898 & {6906}  &  6900 & 6900  \\
G77 & 8124  &9108 &9522  &9810  &{9812}  &9804   & 9804   \\
G81 & 11338  & 12886  & 13836  & 13840  & {13852}  & 13846  & 13846  \\
\bottomrule[0.5pt]
\end{tabular}}
\end{table}

\subsection{Discussion}


Both baseline algorithms, \texttt{CRA} and \texttt{FEM}, implement a form of \textit{continuation methods}~\citep{seguin2022continuation}. They minimize the penalized objective $f(\vx) + \mu \sum_i g(x_i)$ while gradually annealing the penalty parameter~$\mu$.
Although this resembles a \textit{primal-dual} treatment of the \textit{single-constraint} formulation
\begin{equation*}
    \min_{\vx\in[0,1]^n} \; f(\vx) \quad \text{s.t.} \quad \sum_{i=1}^n g(x_i)=0,
\end{equation*}

\texttt{PDBO} differs in two fundamental aspects: (i) \textbf{principled update rule}: Continuation-based methods rely on pre-defined annealing schedule for~$\mu$, whereas \texttt{PDBO}'s dual update follows formally from GDA dynamics; (ii) 
\textbf{richer choices of extensions}: 
\texttt{PDBO} considers the extensions in $S\coloneqq\{f_\vy(\cdot):\vy\in\mathbb R^n\}$, while the continuation methods are restricted to $S'\coloneqq\{f_\mu(\cdot):\mu\in\mathbb R\}$ with $f_\mu(\vx)\coloneqq f(\vx)+\mu\sum\limits_{i=1}^ng(x_i)$.
In the QUBO case, the tightest convex extension in the space $\mathcal{S}$ corresponds to the SDP relaxation~\citep{lemarechal1999semidefinite}, yielding a tighter bound than the eigenvalue relaxation~\citep{nohra2021spectral} resulted from $S'$.

\begin{wraptable}[8]{r}{7.5cm}
\vskip -0.25 in
\caption{Ablation study on Max-Cut}
\vskip 0.08 in
\label{tab:ablation}
\centering
\resizebox{\linewidth}{!}{
\begin{tabular}{@{}lcccccccc@{}}
\toprule
\multirow{2}{*}{Method} & \multicolumn{2}{c}{G72 (n=10k)} & \multicolumn{2}{c}{G77 (n=14k)} & \multicolumn{2}{c}{G81 (n=20k)} \\
\cmidrule(lr){2-3} \cmidrule(lr){4-5} \cmidrule(lr){6-7} 
 & Obj $\uparrow$ & Time $\downarrow$ & Obj $\uparrow$ & Time $\downarrow$ & Obj $\uparrow$ & Time $\downarrow$  \\
\midrule
\texttt{CRA} & 6058 & 53.9 & 8720 & 75.9 & 12450 & 120.4 \\
\texttt{FEM}& 6824 & 2.6 & 9688 & 4.0 & 13684 & 7.5 \\
\texttt{PDBO-s} &6880&1.4&9782&1.8&13786&2.5\\
\texttt{PDBO}   & \bf{6906} & 2.0 & \bf{9812} & 2.1 & \bf{13852} & 3.3 \\
\bottomrule
\end{tabular}}
\end{wraptable}
\textbf{Ablation study.} Table~\ref{tab:ablation} validates the advantage of the primal-dual perspective. First, even the single-constraint variant (denoted as \enquote{\texttt{PDBO-s}}) already outperforms continuation-based baselines (\texttt{CRA}, \texttt{FEM}) across all instances. 
Second, by lifting to $n$ separate constraints—thereby enriching the dual space—the full \texttt{PDBO} framework yields further improvements in solution quality. This progressive gain confirms that both the principled dual update and the expanded constraint formulation contribute to higher solution quality.


\section{Conclusions} \label{sec:conclusions}

In this work, we introduce PDBO, a primal-dual framework for unconstrained binary optimization that reformulates the problem as a minimax optimization and solves it with a tailored GDA algorithm. 
We establish convergence guarantees for PDBO and demonstrate its effectiveness through extensive experiments on public benchmarks.
Our framework opens several promising directions for future research. First, while GDA provides strong empirical performance, exploring more sophisticated minimax solvers, e.g., extragradient methods, could yield better solution quality. Second, the primal-dual perspective could be extended to more complex discrete problems with additional constraints, potentially expanding its applicability to broader classes of combinatorial optimization.
The demonstrated effectiveness of PDBO suggests that the reformulation of binary optimization as a minimax problem is a fruitful approach, warranting further investigation into both algorithmic improvements and theoretical understanding.

\bibliography{ref}
\bibliographystyle{unsrt}

\newpage
\appendix

\section{Proofs}
\label{apdx:proofs}

\subsection{Proof of Lemma~\ref{lem:multi_linear_rounding}}
\begin{proof}

By definition, for any $\vx\in[0,1]^n$, $f(\vx)$ is the convex combination of values in $\left\{F(\bm{\xi}):\bm{\xi}\in\{0,1\}^n\right\}$, with the weights given by $P_\vx$.
Therefore, $f(\vx)\geq\underset{\bm{\xi}\in\{0,1\}^n}\min F(\bm{\xi})$ for any $\vx\in[0,1]^n$.
This means $f$ does not introduce extra minimum on $[0,1]^n\backslash\{0,1\}^n$, and hence the conclusion follows evidently.

\end{proof}

\subsection{Proof of Proposition~\ref{prop:multi_linear_of_polynomial}}
\begin{proof}
\begin{equation*}
\begin{aligned}
    f(\vx)&=\E_{\bm{\xi}\sim P_\vx}\left(\sum\limits_{\vd}a_\vd\xi_1^{d_1}...\xi_n^{d_n}\right)&&\#\text{by definition}\\
    &=\sum\limits_{\vd}a_\vd\left(\E_{\bm{\xi}\sim P_\vx}\xi_1^{d_1}...\xi_n^{d_n}\right)&&\#\;\text{interchange }\E\text{ and }\sum\\
    &=\sum\limits_{\vd}a_\vd\left(\E_{\bm{\xi}\sim P_\vx}\xi_1^{d_1}\right)\cdots\left(\E_{\bm{\xi}\sim P_\vx}\xi_n^{d_n}\right)&&\#\;\xi_i\text{ are independent}\\
    &=\sum\limits_{\vd}a_\vd\left(\E_{\bm{\xi}\sim P_\vx}\xi_1\right)\cdots\left(\E_{\bm{\xi}\sim P_\vx}\xi_n\right)&&\#\;\xi_i^{d_i}=\xi_i\;(\text{for } d_i\neq 0)\text{, since }\xi_i\in\{0,1\}\\
    &=\sum\limits_{\vd}a_\vd\prod\limits_{i:d_i\neq 0}x_i&&\#\;\E_{\bm{\xi}\sim P_\vx}\xi_i=x_i
\end{aligned}
\end{equation*}
\end{proof}

\subsection{Proof of Proposition~\ref{prop:lagrangian_primal}}
\begin{proof}
Observe that
\begin{equation*}
\begin{aligned}
    \underset{\vy\in\mathbb R^n}\sup\;L(\vx,\vy)&=f(\vx)+\underset{\vy\in\mathbb R^n}\sup \;\sum\limits_{i=1}^n\evy_ig(\evx_i)\\
   & =\begin{cases}
    f(\vx),&\text{if } \vx\in\{0,1\}^n\\ +\infty,&\text{otherwise}
\end{cases}
\end{aligned}
\end{equation*}
Therefore, $\underset{\vx\in[0,1]^n}{\inf}\underset{\vy\in\mathbb R^n}{\sup}\;L(\vx,\vy)=\underset{\vx\in\{0,1\}^n}{\inf}\;f(\vx)=v^*$

\end{proof}

\subsection{Proof of Proposition~\ref{prop:y_bounded}}
\begin{proof}
The three conditions for $g$ derive the following properties:
\begin{itemize}
    \item $g(\frac12)\leq g(x)< 0$ for any $x\in(0,1)$
    \item $g'(x)$ is monotone non-decreasing
    \item $g'(\frac12-\delta)<0<g'(\frac12+\delta)$ and $g'(\frac12-\delta)+g'(\frac12+\delta)=0$
   
\end{itemize}

    For any $t \geq 0$, the update $y_i^{t+1} - y_i^t = \beta \cdot g(\evx_i^t)$ lies within $[\beta\cdot g(\frac12), 0]$. Consequently, the sequence $\{y_i^t\}_{t\geq 0}$ is monotone non-increasing.
    If $\evy_i^t> \frac{1+\Theta}{g'(\frac12-\delta)}$ for any $t\geq 0$, there is nothing left to prove, otherwise we denote by $T$ the smallest time index such that $\evy_i^T\leq \frac{1+\Theta}{g'(\frac12-\delta)}$ and $\evx_i^T\notin[\frac12-\delta,\frac12+\delta]$.
    We have $\evy_i^T\in[2\beta\cdot g(\evx_i^t)+\frac{1+\Theta}{g'(\frac12-\delta)},\frac{1+\Theta}{g'(\frac12-\delta)}]$, since Algorithm~\ref{alg:PDBO} prevents the stagnation of $\evx_i^t$ within $[0.5-\delta,0.5+\delta]$ for any two successive steps.
    Now we discuss if $\evx_i^T$ lies in $[0,\frac12-\delta]$ or $[\frac12+\delta,1]$.
    \begin{itemize}
        \item[\textbf{Case I.}] If $\evx_i^{T}\in[0,\frac12-\delta]$, then $g'(\evx_i^T)\leq g'(\frac12-\delta)<0$, and hence
    \begin{equation*}
    \begin{aligned}
        \frac{\partial}{\partial \evx_i} L(\vx^T,\vy^T)&=\frac{\partial}{\partial \evx_i} f(\vx^T)+\evy_i^T\cdot g'(\evx_i^T) \\
        &\geq -\Theta+\frac{1+\Theta}{g'(\frac12-\delta)}\cdot g'(\frac12-\delta)\\
        &=1
    \end{aligned}
    \end{equation*}
    For the primal update
    \begin{equation*}
    \begin{aligned}
        \evx_i^{T+1}&=\Pi_{[0,1]}\left[\evx_i^T-\alpha\cdot\frac{\partial}{\partial \evx_i} L(\vx^T,\vy^T)\right]\\
        &\leq \Pi_{[0,1]}\left[\evx_i^T-\alpha\right]\\
        &\leq \evx_i^T
    \end{aligned}
    \end{equation*}
    Using an analogous argument, we induce that the sequence $\{\evx_i^{T+t}\}_{t\geq 0}$ is non-increasing and decreases at least $\alpha$ in each step until it reaches $0$ after at most $\frac{\frac12-\delta}{\alpha}\leq \lceil\frac1{2\alpha}\rceil$ steps.
    Therefore, we conclude that both $\evx_i^t$ and $\evy_i^t$ remain unchanged for $t\geq T+\lceil\frac1{2\alpha}\rceil$, and 
    \begin{equation*}
    \begin{aligned}
        \evy_i^{T+\lceil \frac1{2\alpha}\rceil}&= \evy_i^T+\sum\limits_{t=1}^{\lceil\frac1{2\alpha}\rceil}\left(\evy_i^{T+t}-\evy_i^{T+t-1}\right)\\
        &\geq \left(2\beta\cdot g(\evx_i^T)+\frac{1+\Theta}{g'(\frac12-\delta)}\right) + \left(\lceil\frac1{2\alpha}\rceil\cdot \beta\cdot g(\frac12)\right)\\
        &\geq \frac{1+\Theta}{g'(\frac12-\delta)} + \left(2+\ceil{\frac1{2\alpha}}\right)\cdot \beta\cdot g(\frac12) &&&\#g(\evx_i^T)\geq g(\frac12)
    \end{aligned}
    \end{equation*}

    \item[\textbf{Case II.}] If $\evx_i^{T}\in[\frac12+\delta, 1]$, then $g'(\evx_i^T)\geq g'(\frac12+\delta)>0$, and hence
    \begin{equation*}
    \begin{aligned}
        \frac{\partial}{\partial \evx_i} L(\vx^t,\vy^t)&=\frac{\partial}{\partial \evx_i} f(\vx^t)+\evy_i^T\cdot g'(\evx_i^T) \\
        &\leq \Theta+\frac{1+\Theta}{g'(\frac12-\delta)}\cdot g'(\frac12+\delta)\\
        &=-1
    \end{aligned}
    \end{equation*}
    For the primal update
    \begin{equation*}
    \begin{aligned}
        \evx_i^{T+1}&=\Pi_{[0,1]}\left[\evx_i^T-\alpha\cdot\frac{\partial}{\partial \evx_i} L(\vx^t,\vy^t)\right]\\
        &\geq \Pi_{[0,1]}\left[\evx_i^T+\alpha\right]\\
        &\geq \evx_i^T
    \end{aligned}
    \end{equation*}
    Using an analogous argument, we induce that the sequence $\{\evx_i^{T+t}\}_{t\geq 0}$ is non-decreasing and increases at least $\alpha$ in each step until it reaches $1$ after at most $\frac{\frac12-\delta}{\alpha}\leq \lceil\frac1{2\alpha}\rceil$ steps.
    Therefore, we conclude that both $\evx_i^t$ and $\evy_i^t$ remain unchanged for $t\geq T+\lceil\frac1{2\alpha}\rceil$, and 
    \begin{equation*}
    \begin{aligned}
        \evy_i^{T+\lceil \frac1{2\alpha}\rceil}&= \evy_i^T+\sum\limits_{t=1}^{\lceil\frac1{2\alpha}\rceil}\left(\evy_i^{T+t}-\evy_i^{T+t-1}\right)\\
        &\geq \left(2\beta\cdot g(\evx_i^T)+\frac{1+\Theta}{g'(\frac12-\delta)}\right) + \left(\lceil\frac1{2\alpha}\rceil\cdot \beta\cdot g(\frac12)\right)\\
        &\geq \frac{1+\Theta}{g'(\frac12-\delta)} + \left(2+\lceil\frac1{2\alpha}\rceil\right)\cdot \beta\cdot g(\frac12)
    \end{aligned}
    \end{equation*}

    \end{itemize}
\end{proof}

\subsection{Proof of Theorem~\ref{thm:strong_duality}}
\begin{proof}
Clearly, we have the following weak duality:
\begin{equation*}
    \underset{\vx\in[0,1]^n}{\inf}\underset{\vy\in\mathbb R^n}{\sup}\;L(\vx,\vy)\geq \underset{\vy\in\mathbb R^n}{\sup}\underset{\vx\in[0,1]^n}{\inf}\;L(\vx,\vy).
\end{equation*}
Taking $\bar \vy\coloneqq \vzero$ yields $\underset{\vx\in[0,1]^n}{\inf}L(\vx,\bar \vy)=\underset{\vx\in[0,1]^n}{\inf} f(\vx)$.
Therefore, we have
\begin{equation*}
\begin{aligned}
\underset{\vy\in\mathbb R^n}{\sup}\underset{\vx\in[0,1]^n}{\inf}\;L(\vx,\vy)&\geq \underset{\vx\in[0,1]^n}{\inf}\;L(\vx,\bar \vy)\\
&=\underset{\vx\in[0,1]^n}{\inf} f(\vx)\\
&=\underset{\vx\in\{0,1\}^n}{\min}f(\vx)&&\left(\text{Lemma~\ref{lem:multi_linear_rounding}}\right)\\
&=\underset{\vx\in[0,1]^n}{\inf}\underset{\vy\in\mathbb R^n}{\sup}\;L(\vx,\vy)&&\left(\text{Proposition~\ref{prop:lagrangian_primal}}\right)
\end{aligned}
\end{equation*}
\end{proof}

\subsection{Proof of Theorem~\ref{thm:convergence_rate}}
\begin{proof}
    Given the update $y_i^{t+1}=y_i^t+\beta\cdot g(\evx_i^t)$, we have $\sum\limits_{t=0}^{T-1}-\beta\cdot g(\evx_i^t)=\evy_i^0-\evy_i^T\leq \evy_i^0-\evy_i^*$ for any $i\in[n]$ and any $T\geq 1$. Consequently,
\begin{align*}
    \frac{\sum\limits_{t=0}^{T-1}\sum\limits_{i=1}^n-g(\evx_i^t)}{T}\leq \frac{||\vy^0-\vy^*||_1 }{\beta\cdot T}  
\end{align*}
Therefore, the number of Algorithm~\ref{alg:PDBO} to return an $\epsilon$-binary point is bounded by $\mathcal{O}\left(\frac{||\vy^0-\vy^*||_1}{\beta\epsilon}\right)$.
\end{proof}

\section{Generalization to the Max-$k$-Cut problem}
\label{apdx:maxkcut}
Given a graph with $n$ nodes, the Max-$k$-Cut problem aims to partition the nodes into $k$ categories that maximize the number of edges connecting the nodes from distinct categories:
\begin{equation*}
\begin{aligned}
    &\underset{\mX\in\mathbb R^{k\times n}}\max && \frac12 \sum\limits_{i=1}^n\sum\limits_{j=1}^n\emW_{i,j}\left(1-\mX_{:,i}^\top \mX_{:,j}\right)\\
    &\quad \text{s.t.}&&\mX_{:,i}\in\left\{\ve^{(1)},\ve^{(2)},...,\ve^{(k)}\right\},\quad\forall i\in[n]
\end{aligned}
\end{equation*}
where $\mW\in\mathbb R^{n\times n}$ is a symmetric matrix representing the weights of edges and each node $i$ is assigned with a one-hot vector $\mX_{:,i}$ representing its category.
The formulation can be equivalently expressed as:
\begin{equation*}
\begin{aligned}
    &\underset{\mX\in\mathbb R^{k\times n}}\min && \text{Tr}\left(\mX\mW\mX^\top\right)\\
    &\quad \text{s.t.}&&\mX_{:,i}\in\left\{\ve^{(1)},\ve^{(2)},...,\ve^{(k)}\right\},\quad\forall i\in[n]
\end{aligned}
\end{equation*}
Note that the domain of each $\mX_{:,i}$ can be equivalently expressed as:
\begin{equation*}
    \left\{\ve^{(1)},\ve^{(2)},...,\ve^{(k)}\right\}=\Delta_k\cap\left\{\mX_{:,i}\in\mathbb R^{k }:g(\mX_{:,i})=0\right\}
\end{equation*}
where $\Delta_k\coloneqq\{\vx\in\mathbb R^k_+:\sum\limits_{i=1}^k\evx_i=1\}$ denotes the $k$-dimensional probability simplex, and $g$ is the convex constraint function as described in Section~\ref{subsec:primal_dual_framework}.
We summarize two possible construction of $g$:
\begin{itemize}
    \item $g(\mX_{:,i})\coloneqq \sum\limits_{j=1}^k \mX_{j,i}^2 - 1$\hfill (quadratic)
    \item $g(\mX_{:,i})\coloneqq \sum\limits_{j=1}^k\mX_{j,i}\log(\mX_{j,i})$\hfill (entropy)
\end{itemize}

We can now reformulate the Max-$k$-Cut problem as:
\begin{equation*}
\begin{aligned}
&\underset{\mX\in\Delta_k^n}{\min}&&\text{Tr}\left(\mX\mW\mX^\top\right)\\
&\quad\text{s.t.}&&g(\mX_{:,i})=0\quad\forall i\in[n]
\end{aligned}
\end{equation*}
The Lagrangian is defined as:
\begin{equation*}
    L(\mX,\vy)\coloneqq \text{Tr}\left(\mX\mW\mX^\top\right)+\sum\limits_{i=1}^ny_ig(\mX_{:,i})
\end{equation*}

To optimize the corresponding minimax problem $\underset{\mX\in\Delta_k^n}{\min}\underset{\vy\in\mathbb R^n}{\max}\;L(\mX,\vy)$, a straightforward approach is to align with Algorithm~\ref{alg:PDBO} by projecting $\mX^t$ onto $\Delta_k^n$.
While there are a variety of methods for normalization, we adopt an alternative by re-parameterizing $\mX_{:,i}\in\Delta_k$ via the softmax function:
\begin{equation*}
\mX_{:,i}=\text{SoftMax}(\mZ_{:,i})=\left(\frac{\exp(\mZ_{j,i})}{\sum\limits_{j=1}^k\exp(\mZ_{j,i})}\right)_{j=1,...,k}
\end{equation*}
and the minimax problem is converted to 
\begin{equation*}
\begin{aligned}
    &\underset{\mZ\in\mathbb R^{k\times n}}{\min}\underset{\vy\in\mathbb R^n}{\max}&&L(\text{SoftMax}(\mZ),\vy)
\end{aligned}
\end{equation*}
Consequently, we can apply the approach described in Section~\ref{subsec:minimax_optimization} to the Max-$k$-Cut problem.

\section{Hyper-parameters}
\label{apdx:hyper_parameters}
\subsection{Experimental Settings}
The multiple initial solutions are selected as follows:
The initial dual variables is set to $\vy^0\coloneqq \bar y\cdot\vone$ and the initial primal variable is sampled from the uniform distribution $\vx^0\sim \text{Uniform}([0,1]^n)$.
The hyper-parameters can be efficiently found via a grid search on small-scale problems.
Table~\ref{tab:hyper_parameters} summarizes the hyper-parameters employed in our experiments.\\\\

\begin{table}[h]
    \centering
    \caption{Hyper-parameters}
    \label{tab:hyper_parameters}
    \begin{tabular}{cccccc} 
    \toprule
    Hyper-parameter & Notation & Max-Cut & MIS & Max-$k$-SAT & Max-$k$-Cut \\
    \midrule
    Number of initial solutions & $B$ &100 &10 &10 &100 \\
    Initial dual variable & $\bar y$ &6  &5  & 2 & 6  \\
    Primal step size & $\alpha$ & 0.025 & 0.02  &0.01  &0.01  \\
    Dual step size &  $\beta$  & 0.025   & 0.02  & 0.005  &0.01  \\
    \bottomrule
    \end{tabular}
\end{table}

\subsection{Sensitivity Analysis}
\label{apdx:sensitivity_analysis}
\begin{table}[h]
\caption{Effect of the initial dual magnitude $\bar y$ on Max-Cut.}
\label{tab:landscape_ablation_apdx}
\centering
\resizebox{\linewidth}{!}{
\begin{tabular}{lccccccc}
\toprule[0.5pt]
$\bar y$ & 0 & 1 & 2 & 4 & 6 & 8 & 10 \\
\midrule[0.4pt]
G70 obj. & 8756  & 9145 ($\times$) & 9398 ($\times$) & 9535 ($\surd$) & 9537 ($\surd$) & 9536 ($\surd$) & 9536 ($\surd$) \\
G81 obj. & 11338 ($\times$) & 12886 ($\times$) & 13836 ($\times$) & 13840 ($\surd$) & 13852 ($\surd$) & 13846 ($\surd$) & 13846 ($\surd$) \\
\bottomrule[0.5pt]
\end{tabular}}
\end{table}
Table~\ref{tab:landscape_ablation_apdx} supports the landscape interpretation through an ablation on the initialization magnitude $\bar y$: performance improves sharply once $\bar y$ is large enough to satisfy the initial convexity condition.

\textbf{Sensitivity to primal step size.}
Table~\ref{tab:sensitivity_alpha} reports the sensitivity to $\alpha$ with $\beta=0.025$ fixed.
\texttt{PDBO} is stable over a broad range of values, except when $\alpha$ is much smaller than $\beta$, in which case the landscape induced by $\vy^t$ evolves too quickly relative to the primal update.

\begin{table}[h]
\caption{Sensitivity to the primal step size $\alpha$ on Max-Cut, with $\beta=0.025$.}
\label{tab:sensitivity_alpha}
\centering
\resizebox{\linewidth}{!}{
\begin{tabular}{lccccccc}
\toprule[0.5pt]
$\alpha$ & 0.0025 & 0.005 & 0.01 & 0.025 & 0.05 & 0.1 & 0.25 \\
\midrule[0.4pt]
G70 obj. & 9365 & 9451 & 9530 & 9537 & 9536 & 9535 & 9539 \\
G81 obj. & 13828 & 13842 & 13848 & 13852 & 13842 & 13854 & 13840 \\
\bottomrule[0.5pt]
\end{tabular}}
\end{table}

\textbf{Sensitivity to dual step size.}
Table~\ref{tab:sensitivity_beta} reports the sensitivity to $\beta$ with $\alpha=0.025$ fixed.
The method remains robust when the dual update is not too aggressive, while overly large $\beta$ degrades solution quality.
This is consistent with the iteration complexity in Theorem~\ref{thm:convergence_rate}, which scales as $\mathcal O(\|\vy^0-\vy^*\|_1/(\beta\epsilon))$, revealing a practical speed-quality trade-off.

\begin{table}[h]
\caption{Sensitivity to the dual step size $\beta$ on Max-Cut, with $\alpha=0.025$.}
\label{tab:sensitivity_beta}
\centering
\resizebox{\linewidth}{!}{
\begin{tabular}{lcccccccc}
\toprule[0.5pt]
$\beta$ & 0.001 & 0.0025 & 0.005 & 0.01 & 0.025 & 0.05 & 0.1 & 0.25 \\
\midrule[0.4pt]
G70 obj. & 9538 & 9536 & 9542 & 9534 & 9537 & 9509 & 9495 & 9415 \\
G81 obj. & 13880 & 13868 & 13862 & 13858 & 13852 & 13840 & 13724 & 13480 \\
\bottomrule[0.5pt]
\end{tabular}}
\end{table}

\textbf{Sensitivity to batch size.}
Table~\ref{tab:sensitivity_batch} reports the sensitivity to the number of parallel initializations $B$.
Larger batches improve objectives at the cost of runtime, while the overall sensitivity remains modest.

\begin{table}[h]
\caption{Sensitivity to batch size $B$ on Max-Cut. Runtime is reported in seconds.}
\label{tab:sensitivity_batch}
\centering
\resizebox{\linewidth}{!}{
\begin{tabular}{lcccc}
\toprule[0.5pt]
$B$ & 1 & 10 & 100 & 1000 \\
\midrule[0.4pt]
G70 obj. (time) & 9525 (0.6) & 9527 (0.7) & 9537 (1.9) & 9540 (9.6) \\
G81 obj. (time) & 13828 (0.9) & 13842 (1.1) & 13852 (3.3) & 13858 (35.8) \\
\bottomrule[0.5pt]
\end{tabular}}
\end{table}

\newpage

\section{Additional Experimental Results}
\label{apdx:additional_experiments}
\subsection{Complete results on Gset instances for Max-Cut and Max-k-cut}
\label{apdx:additional_experiments_mc}
We provide complete results on the Gset instances.
Table~\ref{tab:mc_gset_complete} exhibits the complete results on Gset instances for the Max-Cut problem, while Table~\ref{tab:maxkcut_gset_complete} exhibits the complete results on Gset instances for the Max-3-Cut problem.

Notably, as \texttt{PDBO} exhibits rapid convergence, its solutions can be effectively refined by using them as initial points for \texttt{ABS2}—a method that combines multiple local search strategies and demonstrates strong performance given sufficient computation time. 
This hybrid approach, denoted as \texttt{PDBO+ABS2}, leverages the fast convergence of \texttt{PDBO} with the refinement capability of \texttt{ABS2}.

\begin{sidewaystable}
  \centering
  \caption{Computational results on Gset instances for Max-Cut.}
  \label{tab:mc_gset_complete}
\resizebox{\textwidth}{!}{%
\begin{tabular}{ccccccccccccccccccc}
\toprule
\multirow{2}{*}{Instance} & \multirow{2}{*}{$\vert\mathcal{V}\vert$} & \multirow{2}{*}{$\vert\mathcal{E}\vert$} & \multicolumn{2}{c}{\texttt{PIGNN}} & \multicolumn{2}{c}{\texttt{CRA}} & \multicolumn{2}{c}{\texttt{ROS}} & \multicolumn{2}{c}{\texttt{ANYCSP}} & \multicolumn{2}{c}{\texttt{FEM}} & \multicolumn{2}{c}{\texttt{ABS2}} & \multicolumn{2}{c}{\texttt{PDBO}} & \multicolumn{2}{c}{\texttt{PDBO+ABS2}} 
 \\ \cmidrule(l){4-5}  \cmidrule(l){6-7} \cmidrule(l){8-9} \cmidrule(l){10-11} \cmidrule(l){12-13} \cmidrule(l){14-15} \cmidrule(l){16-17} \cmidrule(l){18-19}
&& & Obj. $\uparrow$    & Time (s) $\downarrow$  & Obj. $\uparrow$    & Time (s) $\downarrow$& Obj. $\uparrow$    & Time (s) $\downarrow$& Obj. $\uparrow$    & Time (s) $\downarrow$& Obj. $\uparrow$    & Time (s) $\downarrow$& Obj. $\uparrow$    & Time (s) $\downarrow$& Obj. $\uparrow$    & Time (s) $\downarrow$& Obj. $\uparrow$    & Time (s) $\downarrow$\\ \midrule  
G1 & 800 & 19176 & 10891& 14.3 & 11398& 25.0 & 11313& 0.8& 11623& 6.5& 11624& 0.8& 11624& 0.1& 11624& 1.6& 11624& 1.6\\ 
G2 & 800 & 19176 & 10930& 23.3 & 11498& 23.7 & 11417& 1.1& 11612& 10.5& 11617& 0.8& 11620& 0.2& 11617& 1.4& 11620& 3.6\\ 
G3 & 800 & 19176 & 10958& 0.2 & 11464& 21.5 & 11383& 0.7& 11620& 3.8& 11622& 0.8& 11622& 0.2& 11621& 1.0& 11622& 1.2\\ 
G4 & 800 & 19176 & 10973& 0.7 & 11486& 24.0 & 11357& 0.6& 11644& 9.8& 11646& 0.7& 11646& 1.1& 11646& 0.9& 11646& 0.9\\ 
G5 & 800 & 19176 & 10901& 0.2 & 11463& 24.8 & 11410& 0.8& 11630& 7.3& 11627& 0.7& 11631& 0.1& 11631& 1.7& 11631& 1.7\\ 
G6 & 800 & 19176 & 1751& 33.3 & 1747& 46.7 & 1949& 0.9& 2173& 13.5& 2178& 0.8& 2178& 0.1& 2178& 1.4& 2178& 1.4\\ 
G7 & 800 & 19176 & 1548& 19.0 & 1512& 48.1 & 1772& 0.9& 1997& 3.7& 2001& 0.8& 2006& 7.1& 2001& 2.2& 2001& 2.2\\ 
G8 & 800 & 19176 & 1588& 58.5 & 1490& 48.1 & 1713& 1.1& 2004& 10.2& 2004& 0.8& 2005& 0.2& 2003& 1.3& 2005& 5.6\\ 
G9 & 800 & 19176 & 1621& 26.8 & 1544& 46.8 & 1784& 0.8& 2044& 13.5& 2051& 0.8& 2054& 0.1& 2054& 1.4& 2054& 1.4\\ 
G10 & 800 & 19176 & 1635& 21.0 & 1481& 48.1 & 1755& 1.3& 1991& 4.8& 1999& 0.8& 2000& 3.0& 1999& 0.9& 1999& 0.9\\ 
G11 & 800 & 1600 & 478& 3.1 & 440& 45.2 & 492& 0.7& 558& 1.4& 558& 0.3& 564& 0.1& 562& 0.6& 564& 1.7\\ 
G12 & 800 & 1600 & 472& 3.6 & 428& 44.7 & 476& 0.4& 548& 0.3& 552& 0.4& 556& 0.5& 554& 0.5& 554& 0.5\\ 
G13 & 800 & 1600 & 490& 4.6 & 430& 44.2 & 514& 0.6& 572& 0.4& 576& 0.3& 582& 0.5& 582& 0.5& 582& 0.5\\ 
G14 & 800 & 4694 & 2815& 12.3 & 2987& 39.2 & 2941& 0.7& 3058& 6.2& 3058& 0.5& 3064& 2.1& 3054& 0.8& 3064& 22.9\\ 
G15 & 800 & 4661 & 2746& 3.8 & 2951& 38.4 & 2944& 0.7& 3040& 2.8& 3047& 0.5& 3045& 1.2& 3041& 0.9& 3043& 2.6\\ 
G16 & 800 & 4672 & 2735& 1.8 & 2966& 39.5 & 2955& 0.8& 3040& 3.4& 3049& 0.5& 3051& 1.2& 3044& 0.9& 3046& 79.9\\ 
G17 & 800 & 4667 & 2699& 40.8 & 2948& 39.2 & 2952& 0.8& 3037& 2.2& 3043& 0.5& 3042& 63.6& 3037& 0.8& 3041& 8.5\\ 
G18 & 800 & 4694 & 833& 14.9 & 822& 44.9 & 859& 1.1& 989& 4.2& 991& 0.5& 991& 33.0& 987& 0.8& 992& 8.1\\ 
G19 & 800 & 4661 & 721& 19.9 & 771& 45.8 & 776& 0.8& 899& 1.5& 905& 0.5& 906& 0.7& 902& 0.7& 902& 0.7\\ 
G20 & 800 & 4672 & 668& 7.7 & 771& 44.3 & 838& 1.2& 936& 1.7& 941& 0.5& 941& 0.1& 940& 1.0& 941& 1.2\\ 
G21 & 800 & 4667 & 753& 11.1 & 766& 46.7 & 810& 0.5& 925& 5.7& 930& 0.4& 931& 1.1& 924& 0.7& 931& 2.5\\ 
G22 & 2000 & 19990 & 12308& 13.1 & 13116& 35.4 & 12920& 1.0& 13347& 13.9& 13356& 0.8& 13357& 10.4& 13347& 2.4& 13359& 4.8\\ 
G23 & 2000 & 19990 & 12288& 2.1 & 13109& 35.0 & 12996& 1.4& 13328& 7.5& 13338& 0.8& 13334& 1.4& 13327& 2.3& 13337& 41.6\\ 
G24 & 2000 & 19990 & 12283& 12.1 & 13067& 35.8 & 12935& 1.1& 13316& 9.9& 13329& 0.8& 13327& 32.9& 13317& 2.4& 13317& 2.4\\ 
G25 & 2000 & 19990 & 12251& 12.6 & 13122& 36.2 & 13008& 0.6& 13323& 12.9& 13334& 0.8& 13329& 55.2& 13321& 2.2& 13325& 18.3\\ 
G26 & 2000 & 19990 & 12247& 14.4 & 13028& 36.4 & 12952& 0.8& 13311& 10.8& 13319& 0.8& 13314& 4.7& 13309& 2.3& 13315& 4.1\\ 
G27 & 2000 & 19990 & 2576& 57.2 & 2749& 46.6 & 2895& 1.0& 3314& 16.2& 3339& 0.8& 3331& 8.0& 3326& 1.1& 3326& 1.1\\ 
G28 & 2000 & 19990 & 2581& 122.1 & 2685& 46.6 & 2972& 0.8& 3275& 18.1& 3291& 0.8& 3283& 1.5& 3296& 2.1& 3296& 2.1\\ 
G29 & 2000 & 19990 & 2730& 50.6 & 2763& 46.6 & 2991& 1.1& 3381& 17.6& 3396& 0.8& 3391& 36.4& 3383& 2.4& 3404& 23.9\\ 
G30 & 2000 & 19990 & 2675& 130.2 & 2748& 46.8 & 2951& 1.0& 3399& 6.2& 3409& 0.8& 3411& 0.9& 3402& 1.9& 3402& 1.9\\ 
G31 & 2000 & 19990 & 2601& 40.4 & 2641& 47.6 & 2890& 0.7& 3296& 16.3& 3302& 0.8& 3288& 1.0& 3302& 2.3& 3306& 5.6\\ 
G32 & 2000 & 4000 & 1170& 5.9 & 1156& 44.7 & 1234& 1.3& 1386& 8.7& 1388& 0.5& 1410& 72.6& 1398& 0.6& 1410& 61.7\\ 
G33 & 2000 & 4000 & 1134& 12.0 & 1122& 43.8 & 1220& 0.7& 1354& 3.3& 1360& 0.4& 1382& 52.7& 1370& 0.6& 1380& 94.5\\ 
G34 & 2000 & 4000 & 1122& 14.4 & 1124& 44.8 & 1204& 0.5& 1358& 7.0& 1358& 0.4& 1384& 9.9& 1376& 0.9& 1384& 127.6\\ 
G35 & 2000 & 11778 & 6968& 4.5 & 7453& 39.4 & 7460& 0.7& 7650& 8.0& 7666& 0.7& 7661& 12.1& 7643& 1.2& 7660& 8.2\\ 
\bottomrule
\end{tabular}%
}
\end{sidewaystable}

\begin{sidewaystable}
  \centering
 \resizebox{\textwidth}{!}{%
\begin{tabular}{ccccccccccccccccccc}
\toprule
\multirow{2}{*}{Instance} & \multirow{2}{*}{$\vert\mathcal{V}\vert$} & \multirow{2}{*}{$\vert\mathcal{E}\vert$} & \multicolumn{2}{c}{\texttt{PIGNN}} & \multicolumn{2}{c}{\texttt{CRA}} & \multicolumn{2}{c}{\texttt{ROS}} & \multicolumn{2}{c}{\texttt{ANYCSP}} & \multicolumn{2}{c}{\texttt{FEM}} & \multicolumn{2}{c}{\texttt{ABS2}} & \multicolumn{2}{c}{\texttt{PDBO}} & \multicolumn{2}{c}{\texttt{PDBO+ABS2}} 
 \\ \cmidrule(l){4-5}  \cmidrule(l){6-7} \cmidrule(l){8-9} \cmidrule(l){10-11} \cmidrule(l){12-13} \cmidrule(l){14-15} \cmidrule(l){16-17} \cmidrule(l){18-19}
&& & Obj. $\uparrow$    & Time (s) $\downarrow$  & Obj. $\uparrow$    & Time (s) $\downarrow$& Obj. $\uparrow$    & Time (s) $\downarrow$& Obj. $\uparrow$    & Time (s) $\downarrow$& Obj. $\uparrow$    & Time (s) $\downarrow$& Obj. $\uparrow$    & Time (s) $\downarrow$& Obj. $\uparrow$    & Time (s) $\downarrow$& Obj. $\uparrow$    & Time (s) $\downarrow$\\ \midrule  
G36 & 2000 & 11766 & 6967& 5.7 & 7435& 38.9 & 7423& 1.4& 7644& 11.7& 7659& 0.7& 7666& 31.4& 7636& 1.0& 7646& 5.6\\ 
G37 & 2000 & 11785 & 7006& 6.7 & 7466& 39.9 & 7432& 1.2& 7654& 8.9& 7670& 0.7& 7667& 56.8& 7657& 0.9& 7676& 54.1\\ 
G38 & 2000 & 11779 & 6993& 1.3 & 7450& 39.1 & 7468& 1.8& 7652& 7.6& 7664& 0.7& 7672& 14.2& 7654& 1.5& 7665& 13.6\\ 
G39 & 2000 & 11778 & 1961& 12.3 & 2115& 46.8 & 2113& 1.2& 2384& 6.8& 2396& 0.6& 2385& 1.4& 2387& 2.0& 2387& 2.0\\ 
G40 & 2000 & 11766 & 1915& 104.4 & 2112& 44.5 & 2149& 1.6& 2376& 13.9& 2394& 0.7& 2389& 3.7& 2365& 1.4& 2396& 83.8\\ 
G41 & 2000 & 11785 & 1976& 21.1 & 2108& 46.0 & 2090& 1.8& 2371& 5.3& 2391& 0.6& 2405& 33.4& 2369& 1.5& 2376& 3.3\\ 
G42 & 2000 & 11779 & 1972& 31.3 & 2203& 44.4 & 2191& 0.6& 2462& 3.8& 2464& 0.7& 2468& 2.4& 2450& 1.3& 2476& 87.8\\ 
G43 & 1000 & 9990 & 6210& 22.9 & 6537& 33.8 & 6427& 0.8& 6660& 7.3& 6660& 0.6& 6660& 0.1& 6660& 1.2& 6660& 1.2\\ 
G44 & 1000 & 9990 & 6160& 1.0 & 6546& 33.9 & 6433& 0.6& 6646& 3.3& 6650& 0.6& 6650& 3.3& 6650& 1.3& 6650& 1.3\\ 
G45 & 1000 & 9990 & 6152& 1.7 & 6520& 35.4 & 6473& 0.8& 6646& 2.0& 6653& 0.6& 6654& 7.3& 6653& 1.1& 6653& 1.1\\ 
G46 & 1000 & 9990 & 6122& 2.3 & 6564& 34.1 & 6463& 0.8& 6639& 7.7& 6646& 0.6& 6649& 2.0& 6649& 1.0& 6649& 1.0\\ 
G47 & 1000 & 9990 & 6174& 1.6 & 6503& 35.0 & 6487& 0.8& 6655& 8.9& 6656& 0.6& 6657& 3.2& 6656& 1.6& 6657& 11.2\\ 
G48 & 3000 & 6000 & 5072& 23.1 & 5852& 41.9 & 5624& 0.9& 6000& 1.5& 6000& 0.3& 6000& 0.2& 6000& 0.4& 6000& 0.4\\ 
G49 & 3000 & 6000 & 4992& 16.0 & 5898& 42.1 & 5552& 0.9& 6000& 0.9& 6000& 0.3& 6000& 0.3& 6000& 0.4& 6000& 0.4\\ 
G50 & 3000 & 6000 & 5040& 16.7 & 5756& 42.7 & 5568& 0.8& 5878& 2.1& 5880& 0.4& 5880& 0.4& 5880& 1.1& 5880& 1.1\\ 
G51 & 1000 & 5909 & 3407& 2.3 & 3732& 39.3 & 3711& 1.0& 3832& 5.8& 3843& 0.5& 3845& 31.3& 3834& 0.9& 3841& 5.5\\ 
G52 & 1000 & 5916 & 3458& 1.1 & 3742& 39.2 & 3734& 1.7& 3837& 7.0& 3845& 0.6& 3846& 70.9& 3834& 0.8& 3841& 4.8\\ 
G53 & 1000 & 5914 & 3526& 9.3 & 3739& 39.3 & 3736& 1.5& 3836& 2.5& 3842& 0.6& 3847& 31.5& 3838& 0.9& 3840& 2.2\\ 
G54 & 1000 & 5916 & 3314& 2.4 & 3743& 39.4 & 3722& 1.4& 3838& 7.2& 3844& 0.5& 3845& 35.2& 3836& 1.1& 3851& 174.6\\ 
G55 & 5000 & 12498 & 9059& 84.3 & 10008& 43.1 & 9797& 1.4& 10248& 22.6& 6418& 0.2& 10254& 93.5& 10240& 1.3& 10258& 120.3\\ 
G56 & 5000 & 12498 & 3127& 106.0 & 3414& 47.1 & 3439& 0.9& 3978& 13.7& 127& 0.2& 3970& 63.0& 3956& 1.5& 3970& 136.7\\ 
G57 & 5000 & 10000 & 2794& 18.5 & 2942& 45.9 & 3036& 0.9& 3398& 18.0& 3408& 0.8& 3474& 150.0& 3460& 1.2& 3478& 157.7\\ 
G58 & 5000 & 29570 & 17435& 7.8 & 18695& 41.8 & 18692& 1.4& 19187& 29.7& 19211& 1.3& 19200& 120.2& 19163& 3.6& 19192& 180.2\\ 
G59 & 5000 & 29570 & 4734& 170.8 & 5580& 48.6 & 5298& 0.9& 5978& 33.1& 6036& 1.4& 6007& 51.8& 5967& 2.7& 6006& 98.2\\ 
G60 & 7000 & 17148 & 12382& 32.4 & 13764& 46.1 & 13413& 1.0& 14108& 27.3& 8746& 0.2& 14138& 171.1& 14111& 2.0& 14111& 2.0\\ 
G61 & 7000 & 17148 & 4518& 58.0 & 4951& 49.4 & 5061& 1.0& 5736& 19.9& 326& 0.2& 5707& 30.7& 5717& 1.8& 5717& 1.8\\ 
G62 & 7000 & 14000 & 3878& 16.0 & 4108& 48.6 & 4268& 1.3& 4758& 13.6& 4762& 1.1& 4848& 176.5& 4808& 1.8& 4828& 175.2\\ 
G63 & 7000 & 41459 & 24237& 106.0 & 26183& 43.5 & 26209& 2.3& 26844& 48.0& 26941& 2.0& 26893& 179.4& 26881& 3.9& 26881& 3.9\\ 
G64 & 7000 & 41459 & 6829& 105.0 & 8114& 50.3 & 7665& 1.3& 8599& 45.9& 8660& 1.9& 8622& 14.9& 8591& 4.5& 8591& 4.5\\ 
G65 & 8000 & 16000 & 4446& 35.2 & 4814& 50.9 & 4898& 1.0& 5414& 26.9& 5420& 1.6& 5520& 148.3& 5482& 1.6& 5504& 167.9\\ 
G66 & 9000 & 18000 & 5142& 20.9 & 5488& 51.7 & 5598& 1.5& 6192& 38.2& 6196& 2.0& 6302& 165.9& 6268& 1.8& 6318& 160.2\\ 
G67 & 10000 & 20000 & 5538& 23.7 & 5948& 53.7 & 6144& 1.5& 6772& 39.9& 6782& 2.4& 6880& 156.3& 6872& 2.5& 6894& 175.8\\ 
G70 & 10000 & 9999 & 8534& 25.2 & 9240& 51.7 & 8872& 1.8& 9379& 35.7& 5120& 0.2& 9510& 175.5& 9537& 1.9& 9537& 1.9\\ 
G72 & 10000 & 20000 & 5588& 44.5 & 6058& 53.9 & 6148& 1.2& 6816& 36.1& 6824& 2.6& 6932& 172.2& 6906& 2.0& 6950& 170.9\\ 
G77 & 14000 & 28000 & 7896& 42.1 & 8720& 75.9 & 8746& 2.2& 9686& 53.5& 9688& 4.0& 9824& 171.1& 9812& 2.1& 9840& 5.1\\ 
G81 & 20000 & 40000 & 11078& 157.6 & 12450& 120.4 & 12320& 5.2& 13670& 73.5& 13684& 7.5& 13850& 177.1& 13852& 3.3& 13860& 4.6\\ 
 \bottomrule
\end{tabular}%
}
\end{sidewaystable}

While PDBO consistently delivers superior performance on Max-3-Cut, the highly specialized ABS2 solver is a strong performer on small-scale Max-Cut instances. On the 66 benchmarks with up to 10,000 nodes, PDBO alone is competitive, matching or beating ABS2 in 16 cases. The most effective strategy, however, is a hybrid approach: using PDBO for rapid initialization followed by ABS2's intensive local search. This hybrid method matches or surpasses standalone ABS2 on 44 of the 66 instances, demonstrating a valuable synergy. 

Moreover, note that ABS2 is not applicable to Max-k-Cut and is confined to QUBO/Max-Cut.
In contrast, PDBO's first-order optimization strategy allows it to address a broader class of combinatorial problems with significant efficiency.

\begin{table}
  \centering
  \caption{Computational results on Gset instances for Max-3-Cut.}
  \label{tab:maxkcut_gset_complete}
\resizebox{0.9\textwidth}{!}{%
\begin{tabular}{ccccccccccc}
\toprule
\multirow{2}{*}{Instance} & \multirow{2}{*}{$\vert\mathcal{V}\vert$} & \multirow{2}{*}{$\vert\mathcal{E}\vert$} & \multicolumn{2}{c}{\texttt{ROS}} & \multicolumn{2}{c}{\texttt{ANYCSP}} & \multicolumn{2}{c}{\texttt{FEM}} &  \multicolumn{2}{c}{\texttt{PDBO}} 
 \\ \cmidrule(l){4-5}  \cmidrule(l){6-7} \cmidrule(l){8-9} \cmidrule(l){10-11} 
&& & Obj. $\uparrow$    & Time (s) $\downarrow$  & Obj. $\uparrow$    & Time (s) $\downarrow$& Obj. $\uparrow$    & Time (s) $\downarrow$& Obj. $\uparrow$    & Time (s) $\downarrow$\\ \midrule  
G1 & 800 & 19176 & 14892& 1.2 & 15121& 18.0 & 15058& 1.0& 15144& 3.1 \\ 
G2 & 800 & 19176 & 14892& 1.3 & 15119& 15.8 & 15070& 0.8& 15159& 3.7 \\ 
G3 & 800 & 19176 & 14932& 1.2 & 15129& 17.0 & 15080& 1.1& 15161& 4.2 \\ 
G4 & 800 & 19176 & 14883& 1.3 & 15141& 15.7 & 15080& 0.9& 15182& 4.5 \\ 
G5 & 800 & 19176 & 14859& 1.4 & 15144& 16.8 & 15087& 1.1& 15185& 4.6 \\ 
G6 & 800 & 19176 & 2341& 1.2 & 1364& 10.5 & 2509& 0.8& 2616& 2.9 \\ 
G7 & 800 & 19176 & 2152& 1.5 & 1131& 0.2 & 2306& 0.8& 2393& 3.9 \\ 
G8 & 800 & 19176 & 2083& 1.6 & 1114& 0.2 & 2321& 1.1& 2413& 3.6 \\ 
G9 & 800 & 19176 & 2214& 1.5 & 1264& 0.2 & 2357& 0.8& 2449& 3.2 \\ 
G10 & 800 & 19176 & 2104& 1.4 & 1110& 0.2 & 2297& 0.9& 2393& 4.1 \\ 
G11 & 800 & 1600 & 607& 0.8 & 655& 2.8 & 645& 0.3& 660& 1.9 \\ 
G12 & 800 & 1600 & 600& 0.8 & 646& 5.6 & 633& 0.3& 653& 1.9 \\ 
G13 & 800 & 1600 & 621& 0.8 & 671& 5.4 & 661& 0.3& 677& 1.7 \\ 
G14 & 800 & 4694 & 3886& 0.9 & 3987& 7.3 & 3960& 0.4& 3998& 2.2 \\ 
G15 & 800 & 4661 & 3855& 1.0 & 3966& 6.5 & 3933& 1.1& 3968& 2.1 \\ 
G16 & 800 & 4672 & 3876& 0.9 & 3966& 7.3 & 3934& 0.4& 3971& 2.2 \\ 
G17 & 800 & 4667 & 3861& 0.8 & 3956& 3.9 & 3929& 1.1& 3964& 2.0 \\ 
G18 & 800 & 4694 & 1080& 1.0 & 1110& 7.0 & 1140& 0.3& 1188& 2.2 \\ 
G19 & 800 & 4661 & 939& 0.9 & 971& 7.8 & 1020& 0.6& 1062& 2.3 \\ 
G20 & 800 & 4672 & 966& 0.9 & 1004& 7.7 & 1048& 0.3& 1106& 2.2 \\ 
G21 & 800 & 4667 & 978& 0.9 & 1010& 5.7 & 1063& 0.4& 1102& 2.0 \\ 
G22 & 2000 & 19990 & 16679& 1.5 & 17068& 16.0 & 16964& 0.8& 17122& 4.3 \\ 
G23 & 2000 & 19990 & 16676& 1.3 & 17084& 25.2 & 16954& 0.8& 17127& 4.3 \\ 
G24 & 2000 & 19990 & 16701& 1.4 & 17063& 22.1 & 16959& 1.0& 17119& 4.8 \\ 
G25 & 2000 & 19990 & 16690& 1.4 & 17072& 24.2 & 16971& 1.0& 17115& 4.2 \\ 
G26 & 2000 & 19990 & 16686& 1.4 & 17056& 22.2 & 16941& 0.9& 17113& 4.1 \\ 
G27 & 2000 & 19990 & 3534& 1.4 & 3129& 24.9 & 3814& 0.8& 3968& 4.2 \\ 
G28 & 2000 & 19990 & 3456& 1.4 & 3030& 13.2 & 3752& 0.8& 3921& 4.1 \\ 
G29 & 2000 & 19990 & 3593& 2.1 & 3225& 9.3 & 3899& 2.4& 4056& 4.9 \\ 
G30 & 2000 & 19990 & 3593& 1.6 & 3231& 14.9 & 3887& 2.1& 4078& 5.6 \\ 
G31 & 2000 & 19990 & 3434& 1.4 & 3067& 18.9 & 3773& 0.8& 3956& 5.8 \\ 
G32 & 2000 & 4000 & 1509& 1.1 & 1602& 11.3 & 1591& 0.3& 1635& 2.3 \\ 
G33 & 2000 & 4000 & 1446& 1.0 & 1572& 11.8 & 1553& 1.1& 1602& 2.3 \\ 
G34 & 2000 & 4000 & 1437& 1.0 & 1551& 11.1 & 1546& 0.4& 1591& 2.3 \\ 
G35 & 2000 & 11778 & 9752& 1.2 & 9981& 16.3 & 9897& 1.9& 9997& 3.4 \\ 
 G36 & 2000 & 11766 & 9744& 1.1 & 9978& 18.2 & 9873& 1.9& 10003& 3.5 \\ 
G37 & 2000 & 11785 & 9766& 1.1 & 9989& 17.0 & 9902& 1.8& 10004& 3.1 \\ 
G38 & 2000 & 11779 & 9795& 1.2 & 9986& 4.8 & 9893& 1.9& 10000& 3.1 \\ 
G39 & 2000 & 11778 & 2533& 1.4 & 2525& 9.7 & 2722& 0.5& 2857& 3.8 \\ 
G40 & 2000 & 11766 & 2461& 1.3 & 2536& 16.3 & 2671& 1.9& 2833& 3.3 \\ 
G41 & 2000 & 11785 & 2524& 1.2 & 2546& 17.8 & 2679& 0.5& 2834& 3.8 \\ 
G42 & 2000 & 11779 & 2582& 1.2 & 2631& 18.0 & 2777& 1.9& 2926& 3.4 \\ 
G43 & 1000 & 9990 & 8312& 1.1 & 8525& 4.3 & 8478& 0.6& 8566& 3.3 \\ 
G44 & 1000 & 9990 & 8349& 1.4 & 8523& 12.8 & 8486& 0.6& 8550& 3.1 \\ 
G45 & 1000 & 9990 & 8354& 1.6 & 8512& 11.9 & 8486& 0.5& 8559& 3.0 \\ 
G46 & 1000 & 9990 & 8346& 4.1 & 8508& 7.2 & 8473& 1.5& 8558& 3.0 \\ 
G47 & 1000 & 9990 & 8338& 2.7 & 8523& 10.7 & 8483& 0.5& 8569& 3.1 \\ 
G48 & 3000 & 6000 & 5920& 2.7 & 6000& 7.7 & 6000& 0.4& 6000& 2.3 \\ 
G49 & 3000 & 6000 & 5914& 2.3 & 6000& 6.9 & 6000& 0.4& 6000& 2.4 \\ 
G50 & 3000 & 6000 & 5924& 1.9 & 5998& 9.7 & 5998& 0.5& 6000& 2.7 \\ 
G51 & 1000 & 5909 & 4897& 1.0 & 5004& 8.3 & 4957& 0.4& 5012& 2.3 \\ 
G52 & 1000 & 5916 & 4910& 1.0 & 5019& 9.1 & 4970& 0.4& 5012& 2.3 \\ 
G53 & 1000 & 5914 & 4918& 1.0 & 5007& 4.7 & 4964& 0.4& 5015& 2.2 \\ 
G54 & 1000 & 5916 & 4886& 0.9 & 5015& 8.7 & 4966& 0.4& 5013& 2.2 \\ 
G55 & 5000 & 12498 & 11992& 1.4 & 12351& 27.3 & 12210& 2.2& 12345& 4.3 \\ 
G56 & 5000 & 12498 & 4126& 1.7 & 4495& 31.3 & 4425& 2.2& 4663& 4.0 \\ 
G57 & 5000 & 10000 & 3688& 1.5 & 3940& 29.4 & 3900& 0.4& 4028& 3.6 \\ 
G58 & 5000 & 29570 & 24482& 4.2 & 25027& 46.5 & 24807& 3.7& 25088& 5.4 \\ 
G59 & 5000 & 29570 & 6549& 5.2 & 6316& 43.3 & 6785& 3.9& 7190& 8.0 \\ 
G60 & 7000 & 17148 & 16496& 1.9 & 16978& 41.6 & 16797& 2.9& 16961& 4.5 \\ 
G61 & 7000 & 17148 & 5989& 2.0 & 6464& 41.2 & 6401& 2.9& 6730& 4.8 \\ 
G62 & 7000 & 14000 & 5149& 1.9 & 5487& 37.0 & 5441& 2.5& 5635& 4.3 \\ 
G63 & 7000 & 41459 & 34331& 3.0 & 35087& 49.3 & 34776& 5.3& 35194& 8.0 \\ 
G64 & 7000 & 41459 & 9399& 2.8 & 9219& 48.1 & 9767& 5.1& 10333& 8.4 \\ 
G65 & 8000 & 16000 & 5914& 2.2 & 6255& 44.3 & 6207& 2.9& 6425& 4.5 \\ 
G66 & 9000 & 18000 & 6739& 2.5 & 7147& 39.4 & 7103& 3.2& 7338& 4.5 \\ 
G67 & 10000 & 20000 & 7364& 3.0 & 7797& 55.9 & 7748& 3.3& 8015& 6.0 \\ 
G70 & 10000 & 9999 & 9983& 1.9 & 9909& 16.2 & 9999& 1.4& 9999& 3.3 \\ 
G72 & 10000 & 20000 & 7435& 2.9 & 7906& 58.9 & 7835& 0.7& 8111& 4.7 \\ 
G77 & 14000 & 28000 & 10559& 5.3 & 11158& 84.0 & 11102& 4.4& 11467& 5.0 \\ 
G81 & 20000 & 40000 & 14907& 9.7 & 15727& 115.0 & 15683& 6.1& 16191& 7.9 \\ 

 \bottomrule
\end{tabular}%
}
\end{table}

\newpage

\subsection{Experiments for GPU efficiency}
In this section we conduct comparison between the CPU/GPU versions of PDBO.

\begin{figure}[h]
    \centering
\includegraphics[width=0.5\linewidth]{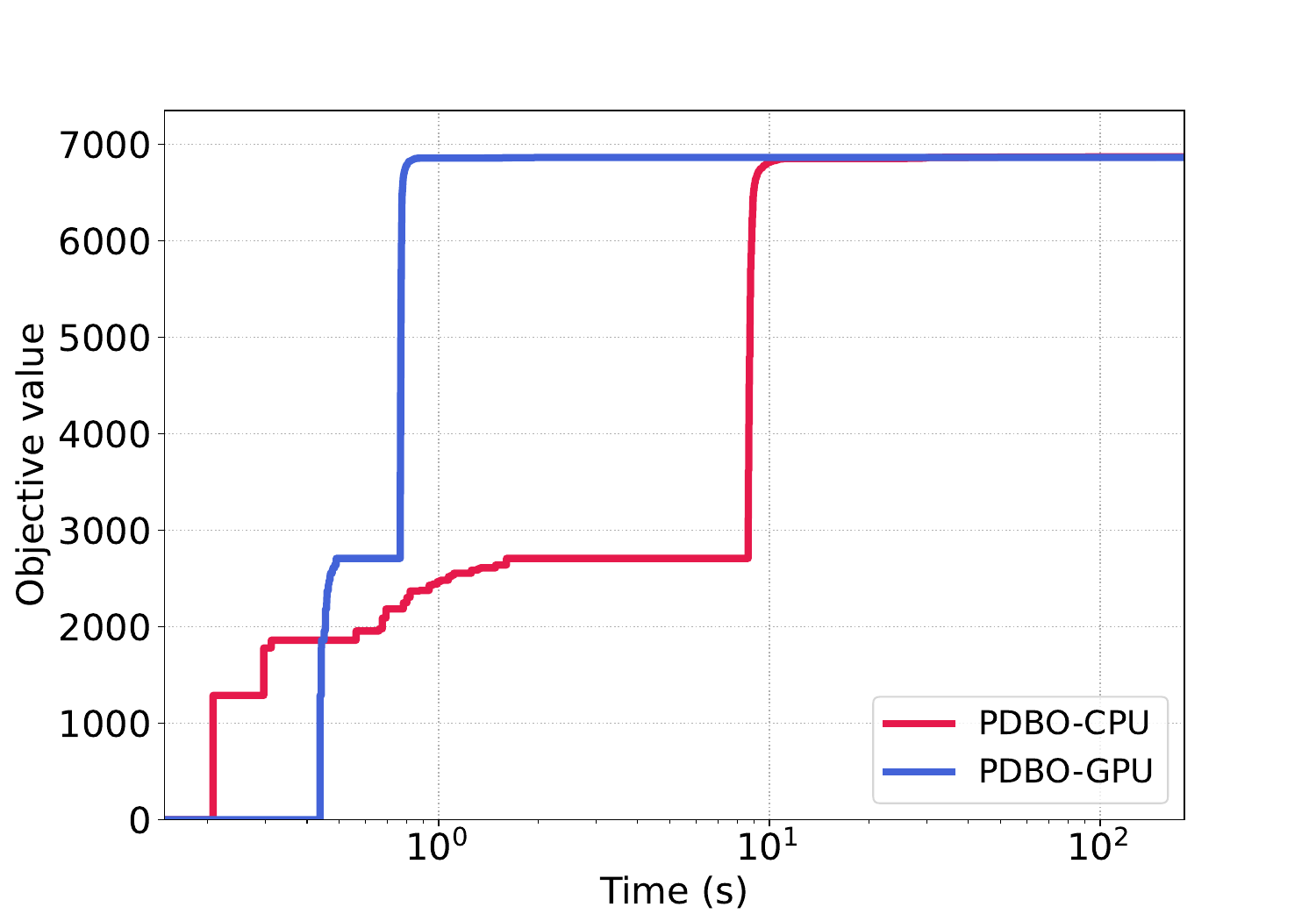}
    \caption{The objective value of Max-Cut on G67, as a function of runtime}
    \label{fig:GPUvsCPU}
\end{figure}

Figure~\ref{fig:GPUvsCPU} compares the solution quality of each method throughout the optimization process, plotting the objective value against runtime under a 180-second time limit, on the G67 instance for Max-Cut.
Clearly, \texttt{PDBO-GPU} (blue) performs much faster than \texttt{PDBO-CPU} (red). The initial latency of \texttt{PDBO-GPU} could be attributed to the \textit{just-in-time} nature of Jax.

\subsection{Experiments on MIS}
Given a graph $G \coloneqq (V, E)$, the MIS problem seeks the largest vertex subset $S \subseteq V$ such that no two vertices in $S$ are adjacent. This can be formulated as:
\begin{equation*}
\begin{aligned}
&\max_{\vx \in \{0,1\}^n} && \sum_{i=1}^n x_i \\
&\quad \; \text{s.t.} && x_i x_j = 0, \quad \forall (i,j) \in E
\end{aligned}
\end{equation*}
where $x_i \in \{0,1\}$ indicates whether vertex $i$ is included in the independent set.
Following the approach of~\cite{schuetz2022combinatorial,ichikawa2024controlling}, we reformulate MIS as a QUBO by penalizing the adjacency constraints:
\begin{equation*}
    \max_{\vx \in \{0,1\}^n} \sum_{i=1}^n x_i - \lambda \cdot \sum_{(i,j) \in E} x_i x_j
\end{equation*}
where $\lambda > 0$ is a penalty parameter. In line with established practice, we set $\lambda=4$ for all methods to ensure a consistent comparison.

For MIS, we employ \textit{d}-regular random graphs (\textit{d}-RRGs), where each node has degree $d$, following the experimental setup of~\cite{schuetz2022combinatorial}. To evaluate performance across both sparse and dense graph structures, we consider degrees $d \in \{3, 100\}$ and graph sizes $n \in \{10^4, 5 \times 10^4\}$ nodes. For each $(n, d)$ configuration, we generate 20 independent graph instances to ensure statistical significance.

The objective of the MIS problem is to find the largest set of non-adjacent nodes. Although we formulate MIS as a QUBO via constraint penalization, we only report the objective value (independent set size) for solutions that satisfy the independence constraints. We evaluate \texttt{PDBO} against baseline algorithms \texttt{PIGNN}, \texttt{CRA}, \texttt{FEM}, \texttt{ABS2} and \texttt{Gurobi} on \textit{d}-RRG instances, with performance averaged over 20 instances per configuration, as shown in Table~\ref{tab:mis}.
The results show that \texttt{PIGNN}, \texttt{CRA}, \texttt{FEM} and even \texttt{Gurobi} struggle on dense graphs with degree 100 (entries marked with \enquote{--} indicate that \texttt{FEM} failed to return feasible solutions). This aligns with the known computational hardness of dense graphs~\citep{barbier2013hard}. In contrast, both \texttt{PDBO} and \texttt{ABS2} achieve strong performance even on dense instances. Notably, \texttt{PDBO} outperforms all baselines, identifying better solutions in less time on both sparse and dense graphs. 

\begin{table}[h]
\caption{Results on $d$-RRGs for MIS}
\label{tab:mis}
\centering
\resizebox{\textwidth}{!}{
\renewcommand{\arraystretch}{1.2}
\begin{tabular}{@{}lcccccccccc@{}}
\toprule
 \multirow{2}{*}{Method}& \multicolumn{2}{c}{$(10^4,3)$} & \multicolumn{2}{c}{$(10^4,100)$} & \multicolumn{2}{c}{$(5\times10^4,3)$} & \multicolumn{2}{c}{$(5\times10^4,100)$} \\
\cmidrule(lr){2-3} \cmidrule(lr){4-5} \cmidrule(lr){6-7} \cmidrule(lr){8-9} 
 & Obj $\uparrow$ & Time $\downarrow$ & Obj $\uparrow$ & Time $\downarrow$ & Obj $\uparrow$ & Time $\downarrow$ & Obj $\uparrow$ & Time $\downarrow$\\
\midrule
\texttt{PIGNN} & 4176.1 $_{\pm 19.3 }$ & 66.7 $_{\pm 14.2 }$ & 7.5 $_{\pm 6.7 }$ & 0.2 $_{\pm 0.1 }$ & 20705.7 $_{\pm 47.1 }$ & 178.9 $_{\pm 1.0 }$ & 15.8 $_{\pm 13.7 }$ & 0.6 $_{\pm 0.5 }$\\
\texttt{CRA} & 4361.4 $_{\pm 8.1 }$ & 52.2 $_{\pm 0.1 }$ & 361.9 $_{\pm 187.6 }$ & 46.6 $_{\pm 23.2 }$ & 166.0 $_{\pm 75.7 }$ & 20.4 $_{\pm 4.0 }$ & 11.1 $_{\pm 9.8 }$ & 0.5 $_{\pm 0.5 }$\\
\texttt{FEM} & 4417.2 $_{\pm 4.4 }$ & 0.7 $_{\pm 0.0 }$ & {--}  & 180.0 $_{\pm 0.0 }$ & 22053.8 $_{\pm 11.1 }$ & 1.3 $_{\pm 0.1 }$ & --  & 180.0 $_{\pm 0.0 }$\\
\texttt{ABS2} & 4394.8 $_{\pm 6.5 }$ & 154.2 $_{\pm 27.6 }$ & 587.4 $_{\pm 3.5 }$ & 140.0 $_{\pm 36.2 }$ & 21307.7 $_{\pm 43.6 }$ & 173.1 $_{\pm 5.7 }$ & 2971.5 $_{\pm 12.8 }$ & 175.0 $_{\pm 3.9 }$\\
\texttt{Gurobi}&{4141.8 $_{\pm  19.4 }$ }&{147.2 $_{\pm 24.8 }$ }& {448.0 $_{\pm  5.1 }$ }&{0.1 $_{\pm  0.0 }$ }&{20144.3 $_{\pm  181.7 }$}& {96.6 $_{\pm  49.9 }$ }& {2235.6 $_{\pm  12.7 }$ }& {0.1 $_{\pm  0.0 }$} \\
\texttt{PDBO} & \textbf{4431.9 $_{\pm 3.1 }$} & 0.6 $_{\pm 0.0 }$ & \textbf{603.0 $_{\pm 3.3 }$} & 3.6 $_{\pm 0.3 }$ & \textbf{22127.3 $_{\pm 10.3 }$} & 1.1 $_{\pm 0.1 }$ & \textbf{3006.2 $_{\pm 6.2 }$} & 14.6 $_{\pm 0.6 }$\\
\bottomrule
\end{tabular}}
\end{table}

\subsection{Experiments on MAP}
\label{sec:map}
In this section we conduct additional experiments on the Maximum a Posterior (MAP) inference in discrete Markov random fields (MRF). 
Let $\mathcal G$ be a graph of $n$ nodes and $\mathcal C$ is the set of cliques, and let $\mathcal S=\mathcal S_1\times...\times\mathcal S_n$ be the assignment domain of the $n$ nodes (i.e. variables). Consider an MRF representing a joint distribution 
\begin{equation*}
    p(s)=\frac1Z\prod\limits_{C\in\mathcal C}\psi_C(s_C), \quad\forall s\in\mathcal S
\end{equation*}
where $s_C$ is the joint assignment of the variables in clique $C$, $\psi_C$ are positive functions called potentials, and $Z$ is the partition function.

The MAP inference problem is to find the most likely assignment to the variables, which can be equivalently modeled as:
\begin{equation*}
    \underset{s\in\mathcal S}{\min} \;\sum\limits_{C\in\mathcal C}-\log \psi_C(s_C)
\end{equation*}
As shown by~\cite{le2018continuous}, when $|\mathcal S_i|=2$ for all $i$, the MAP inference can be reformulated as an binary optimization with polynomial objective, which is applicable to PDBO. Furthermore, one can easily handles the cases where $|\mathcal S_i|>2$ with similar techniques as we discussed in Appendix~\ref{apdx:maxkcut}.

\subsubsection{Protein Prediction}
In particular, we adopt the protein-prediction dataset~\citep{jaimovich2006towards}, which is widely adopted in the related literature and competitions~\citep{kappes2015comparative, elidan2011probabilistic}.
The dataset includes 8 instances and Table~\ref{tab:map} presents the average objective value and running time for each method.
Please refer to~\cite{kappes2015comparative} for detailed description on the baselines as well as the formal problem definition.

\begin{table}[h]
\caption{Results for MAP on the protein-prediction dataset}
\label{tab:map}
    \centering

    \begin{tabular}{ccc}
    \toprule
         & Obj.$\downarrow$& Time$\downarrow$ \\
\midrule
\texttt{ogm-LBP-LF2} &52942.95&69.86\\
\texttt{ogm-LF-3}&57944.06&25.99\\
\texttt{ogm-LBP-0.5}&53798.89&60.97\\
\texttt{ogm-TRBP-0.5}&61386.17&86.03\\
\texttt{ADDD}&106216.86&10.97\\
\texttt{MPLP}&101531.75&86.6\\
\texttt{ogm-LP-LP}&102918.41&180.72\\
\texttt{ogm-ILP}&60047.02&2262.83\\
\texttt{BRAOBB-1}&61079.07&3600.16\\
\texttt{PDBO}&\textbf{52842.61}&15.07\\ 
\bottomrule
    \end{tabular}
    
\end{table}
The results show that \texttt{PDBO} is both efficient and effective for handling the MAP problems.

\subsubsection{UAI2022 Competition}
We also evaluate PDBO on the datasets from the UAI 2022 competition (\url{https://uaicompetition.github.io/uci-2022/}), a standard modern benchmark. 
In particular, we consider the \textbf{Grids} dataset (7 instances in total) that aligns with the binary optimization formulation~(\ref{eq:binary_optimization}).
We compare against three state-of-the-art baselines—\texttt{LBP}, \texttt{TRBP}, and \texttt{Toulbar2}—with their results taken directly from the recent and comprehensive study of~\cite{wang2024neurolifting}:
\begin{itemize}
    \item Loopy Belief Propagation (denoted as \texttt{LBP}) \& Tree-reweighted Belief Propagation (denoted as \texttt{TRBP})~\citep{wainwright2005map}: Classic and widely-used message-passing algorithms.
    \item \texttt{Toulbar2}~\citep{de2023toulbar2}: the winner on all MPE and
MMAP task categories of UAI 2022 Inference Competition.
\end{itemize}

\begin{table}[h]
\caption{Results on the UAI 2022 benchmarks, lower values imply better performance}
    \label{tab:uai2022}
    \centering

    \resizebox{\textwidth}{!}{
    \begin{tabular}{cccccccc}
    \toprule
    & Grids\_{19}
    &Grids\_{21} &Grids\_{24}
    &Grids\_{25}
    &Grids\_{26} &Grids\_{27} &Grids\_{30} \\
    \midrule
     \texttt{LBP} &-2250.440&-13119.300&-13210.400&-2170.890&-2063.350&-9024.640&-2142.890\\
\texttt{TRBP} &-2103.610&-12523.300&-13260.900&-2171.050&-1903.910&-9019.470&-2154.910\\
\texttt{Toulbar2}&-2643.107&-18895.393&-18274.302&-2620.268&-3010.719&-12284.284&-2984.248\\
\texttt{PDBO} & \bf-2696.341&\bf-19235.873&\bf-18691.510&\bf-2653.470&\bf-3021.883&\bf-12481.272&\bf-2989.997\\
\bottomrule
    \end{tabular}}
\end{table}
The results are shown in Table~\ref{tab:uai2022}. \texttt{PDBO} achieves the best performance on all seven instances. Notably, it finds these superior solutions within 5.9 seconds on average, while the baselines from~\cite{wang2024neurolifting} (including \texttt{Toulbar2} with a 1,200-second limit) could not match this performance. These results, conducted on a standard modern benchmark and compared against baseline results, robustly demonstrate PDBO's effectiveness and efficiency.

\subsection{Dual Certificate}
When the initial value $\vy^0$ is sufficiently large (e.g. ${\bar{y}} \geq -\lambda_{min}(W)$ for max-cut), the Lagrangian extension $f_\vy(\vx) \coloneqq L(\vx,\vy)$ turns out to be convex in $\vx$. In this regime, we are able to globally minimize $f_\vy(\vx)$ via classic gradient descent (albeit PDBO employs a single gradient descent step instead), the corresponding optimal value would provide a valid dual bound on the optimal objective value, serving as a dual certificate. We present the corresponding dual bounds of the Max-Cut problem implied by PDBO in Table~\ref{tab:dual_bounds}.
\begin{table}[h]
\caption{Dual bounds for Max-Cut, implied by PDBO}
\label{tab:dual_bounds}
    \centering
    
    \begin{tabular}{cccccc}
    \toprule
    Max-Cut& G67&G70&G72&G77&G81\\
    \midrule
      primal bound   & 6872 &9537&6906&9812&13852 \\
     dual bound    & 8843.7&	13891.5&	8906.1&	12767.7&	18172.2\\
     \bottomrule
    \end{tabular}
\end{table}

\end{document}

%% file: math_commands.tex
\usepackage{amsmath,amsfonts,bm}

\def\vzero{{\bm{0}}}
\def\vone{{\bm{1}}}

\def\vd{{\bm{d}}}
\def\ve{{\bm{e}}}

\def\vx{{\bm{x}}}
\def\vy{{\bm{y}}}

\def\evx{{x}}
\def\evy{{y}}

\def\mI{{\bm{I}}}

\def\mW{{\bm{W}}}
\def\mX{{\bm{X}}}

\def\mZ{{\bm{Z}}}

\DeclareMathAlphabet{\mathsfit}{\encodingdefault}{\sfdefault}{m}{sl}
\SetMathAlphabet{\mathsfit}{bold}{\encodingdefault}{\sfdefault}{bx}{n}

\def\emW{{W}}

\newcommand{\E}{\mathbb{E}}

